\documentclass[a4paper,10pt]{article}

\usepackage{latexsym}
\usepackage{mathrsfs}
\usepackage{amsfonts,amsmath,amssymb}
\usepackage{indentfirst}
\usepackage{subeqnarray}
\usepackage{verbatim}
\usepackage[pdftex]{graphicx}
\usepackage{epstopdf}
\usepackage{stmaryrd}
\usepackage{multirow}
\usepackage{diagbox}
\usepackage{subcaption}

\newcommand{\bx}{\mbox{\boldmath{$x$}}}

\newcommand{\bz}{\mbox{\boldmath{$z$}}}

\newcommand{\btau}{\mbox{\boldmath{$\tau$}}}
\newcommand{\bphi}{\mbox{\boldmath{$\phi$}}}
\newcommand{\bzero}{\mbox{\boldmath{$0$}}}

\newcommand{\bb}{\mbox{\boldmath{$b$}}}

\newcommand{\fb}{\mbox{\boldmath{$f$}}}
\newcommand{\bg}{\mbox{\boldmath{$g$}}}
\newcommand{\bK}{\mbox{\boldmath{$K$}}}

\newcommand{\bn}{\mbox{\boldmath{$n$}}}
\newcommand{\bt}{\mbox{\boldmath{$t$}}}
\newcommand{\bu}{\mbox{\boldmath{$u$}}}

\newcommand{\bI}{\mbox{\boldmath{$I$}}}
\newcommand{\bv}{\mbox{\boldmath{$v$}}}
\newcommand{\bV}{\mbox{\boldmath{$V$}}}

\newcommand{\bW}{\mbox{\boldmath{$W$}}}

\newcommand{\by}{\mbox{\boldmath{$y$}}}
\newcommand{\bcW}{\mathbfcal{W}}
\newcommand{\bcB}{\mathbfcal{B}}
\newcommand{\cU}{\mathcal{U}}
\newcommand{\cE}{\mathcal{E}}
\newcommand{\cN}{\mathcal{N}}
\newcommand{\cR}{\mathcal{R}}
\newcommand{\cT}{\mathcal{T}}
\newcommand{\dS}{\mathbb{S}}
\newcommand{\dQ}{\mathbb{Q}}

\newcommand{\bkappa}{\mbox{\boldmath{$\kappa$}}}
\newcommand{\bvarepsilon}{\mbox{\boldmath{$\varepsilon$}}}
\newcommand{\bsigma}{\mbox{\boldmath{$\sigma$}}}
\newcommand{\bSigma}{\mbox{\boldmath{$\Sigma$}}}

\newcommand{\dx}{\mathrm{d}\bx}
\newcommand{\ds}{\mathrm{d}s}

\newcommand{\Real}{\mathbb{R}}
\newcommand{\set}[1]{\left\{#1\right\}}

\newtheorem{theorem}{Theorem}[section]

\newtheorem{example}{Example}

\newtheorem{remark}[theorem]{Remark}
\numberwithin{equation}{section}

\DeclareMathAlphabet\mathbfcal{OMS}{cmsy}{b}{n}

\usepackage[pdftex,unicode,colorlinks,linkcolor=blue,citecolor=red]{hyperref}

\begin{document}

\begin{center}
\Large\bf Local Randomized Neural Networks with Hybridized Discontinuous Petrov-Galerkin Methods for Stokes-Darcy Flows
\end{center}

	\begin{center}
		{\large\sc Haoning Dang}\footnote{School of Mathematics and Statistics, Xi'an Jiaotong University, Xi'an,
			Shaanxi 710049, P.R. China. E-mail: {\tt haoningdang.xjtu@stu.xjtu.edu.cn}},\quad
		{\large\sc Fei Wang}\footnote{School of Mathematics and Statistics, Xi'an Jiaotong University,
			Xi'an, Shaanxi 710049, China. The work of this author was partially supported by
			the National Natural Science Foundation of China (Grant No.\ 12171383) and 
			Shaanxi Fundamental Science Research Project for Mathematics and Physics (Grant No.\ 22JSY027). Email: {\tt feiwang.xjtu@xjtu.edu.cn}}
	\end{center}

\medskip
\begin{quote}
{\bf Abstract. } 
This paper introduces a new numerical approach that integrates local randomized neural networks (LRNNs) and the hybridized discontinuous Petrov-Galerkin (HDPG) method for solving coupled fluid flow problems. The proposed method partitions the domain of interest into several subdomains and constructs an LRNN on each subdomain. Then, the HDPG scheme is used to couple the LRNNs to approximate the unknown functions. We develop LRNN-HDPG methods based on velocity-stress formulation to solve two types of problems: Stokes-Darcy problems and Brinkman equations, which model the flow in porous media and free flow. We devise a simple and effective way to deal with the interface conditions in the Stokes-Darcy problems without adding extra terms to the numerical scheme. We conduct extensive numerical experiments to demonstrate the stability, efficiency, and robustness of the proposed method. The numerical results show that the LRNN-HDPG method can achieve high accuracy with a small number of degrees of freedom.

\end{quote}

\medskip


\section{Introduction}

The Stokes-Darcy equations describe the flow in a porous medium coupled with a free flow, which has many applications (\cite{Discacciati2009Navier}). To couple the two regions of the Stokes-Darcy flow, the Beavers-Joseph law is a general interface condition derived from experiments (\cite{Beavers1967BoundaryCA}). The Beavers-Joseph-Saffman law is a simplified version of the Beavers-Joseph law (\cite{Saffman1971Boundary}), which leads to a simpler numerical implementation. One numerical challenge for this model is to design a numerical scheme that is stable for both the Stokes-dominated and the Darcy-dominated regions, since the usual Darcy stable numerical methods may diverge for Stokes flow and vice versa. Another challenge is to handle the interface conditions and couple the different numerical methods for different domains,  if we use different numerical methods for each domain to solve the Stokes-Darcy equations. Various traditional methods have been used to tackle these challenges. For instance, \cite{Riviere2005Locally} applied a mixed finite element method in the Darcy region and a discontinuous Galerkin method in the Stokes region. \cite{Riviere2005Analysis} used the interior penalty discontinuous Galerkin method to solve the Stokes-Darcy equations in both regions. \cite{Chen2016Weak} employed the weak Galerkin method for both regions. Another approach is to decouple the Stokes-Darcy equations and solve them separately using different numerical methods. \cite{Mu2007Twogrid} adopted a two-grid method to achieve this decoupling and applied local solvers for each region.

Deep neural networks (DNNs) have attracted a lot of attention since their success in the ILSVRC-2012 competition (\cite{Krizhevsky2012ImageNet}). DNNs have demonstrated remarkable and surprising performance in various machine learning tasks, such as natural language processing, computer vision, electronic design automation, and more. Mathematically, numerical solutions of partial differential equations (PDEs) can be viewed as finding approximations of functions in some function spaces. DNNs are new classes of functions, with powerful approximation properties. 
The universal approximation theorem (\cite{Barron1993Universal,Cybenko1989Approximation}) is a well-known result for DNNs, stating that a single-hidden-layer neural network can approximate any continuous function with a suitable activation function. Furthermore, the analysis in \cite{Barron1994Approximation} shows that the convergence order is independent of the problem dimension. Other studies have explored the approximation capabilities of shallow neural networks (\cite{Xu2020FiniteNeuron,Siegel2020Approximation,Siegel2021SharpLB,E2022Representation,Ming2023BarronClass}). Approximation results for deep neural networks can be found in \cite{E2020MutilLayer} and the references therein.

Neural networks have been applied to solve PDEs in various ways in recent years. One of the most popular approaches is physics-informed neural networks (PINNs, \cite{Raissi2019Physics}), which use a fully connected feedforward neural network and construct the loss functions based on the strong form of PDEs with initial and boundary conditions. Another approach is the deep Galerkin method (\cite{Sirignano2018DGM}), which uses a least-squares loss function. However, the weak form can handle problems with less smooth solutions better than the strong form. Therefore, some methods combine the weak form and neural networks. For example, the deep Ritz method (DRM, \cite{E2018DRM}) uses the energy functional of the weak form as the loss function. Deep Nitsche method (\cite{Liao2021DNM}) uses Nitsche's energy formulation to impose the essential boundary conditions in a different way from DRM. Weak adversarial networks (\cite{Zang2020WAN}) use two networks, one for the trial function and one for the test function, to train the neural networks. All these methods use stochastic gradient descent-type algorithms to solve non-convex optimization problems, which are computationally expensive and time-consuming. Moreover, they may get stuck in local optima or saddle points, resulting in low accuracy and large optimization errors. Hence, these neural network methods have the potential to avoid the curse of dimensionality and solve high-dimensional PDEs, but they are not competitive with traditional methods for low-dimensional problems.

One way to improve the efficiency and reduce the optimization error is to employ randomized neural networks (RNNs, \cite{Igelnik1995Stochastic}), where all parameters except those of the output layer are randomly generated and fixed during the training. The remaining parameters are then obtained by solving a least-squares problem. A special case of RNNs is an extreme learning machine (ELM, \cite{Huang2006Extreme}). It has been shown that a randomized neural network with one hidden layer can approximate any continuous function on a compact domain, provided that a suitable activation function is used and the hidden layer size is sufficiently large (\cite{Igelnik1995Stochastic}). In \cite{Dong2021Local}, local ELM and domain decomposition methods have been combined to solve PDEs, and numerical examples have demonstrated that RNNs can achieve high accuracy with fewer degrees of freedom. Various PDEs have been solved by using ELM-based methods, such as advection-diffusion equations (\cite{Quan2023Solving}), fractional partial differential equations (\cite{Sivalingam2023Novel}), inverse parametric PDEs (\cite{Dong2023Method}), time-dependent PDEs (\cite{Chen2023RFM}), interface problems (\cite{Li2023LRNN}), and more. 
RNNs can also work with the weak form of PDEs and produce good results. For instance, a Petrov-Galerkin method combined with RNNs was proposed in \cite{Shang2023Randomized} to solve linear and nonlinear PDEs, and then the method was applied to solve linear elasticity problems and Navier-Stokes equations in \cite{Shang2023RNNPG}. Local RNNs (LRNNs) with discontinuous Galerkin methods were developed in \cite{Sun2022LRNNDG} for solving PDEs and were extended to solve diffusive-viscous wave equations in \cite{Sun2023DVWE}. Numerical experiments have indicated that these RNNs-based numerical methods have several benefits. They can achieve high accuracy with fewer degrees of freedom, they can naturally and efficiently solve linear and nonlinear time-dependent problems by using the space-time approach, and they can easily handle boundary conditions.

In this paper, we propose a new approach to solve the Stokes-Darcy and Brinkman equations, based on the local randomized neural networks with hybridized discontinuous Petrov-Galerkin (LRNN-HDPG) method. Unlike the conventional discontinuous Petrov-Galerkin method (\cite{Demkowicz2010class, Astaneh2018High}), which uses piecewise discontinuous polynomials for both the trial and the test functions, our method approximates the trial functions by LRNNs, while the test functions are piecewise polynomials. We also present the LRNN-HDG method, which uses LRNNs for both trial and test functions. The proposed method can address the two main challenges of solving the Stokes-Darcy equations in a simple way. First, it is stable for both the Stokes-dominated and the Darcy-dominated regions. Second, it can directly enforce the interface conditions at sampling points on the interface, instead of incorporating them into the numerical scheme. This avoids the difficulty of coupling different numerical methods in different domains. We solve the resulting linear system by the least-squares method. Compared with traditional numerical methods, our method can achieve very accurate numerical solutions with fewer degrees of freedom.

This paper is organized as follows: In Section \ref{Sec:RNN}, we introduce the mathematical representation of RNNs. In Section \ref{Sec:lrnn-hdpg}, we describe the LRNN-HDPG method for the Darcy flows. In Section \ref{Sec:Stokes-Darcy-Brinkman}, we develop the LRNN-HDPG method for the Stokes equations, and then we address the coupled Stokes-Darcy and Brinkman equations. In Section \ref{Sec:NumerExper}, we report some numerical experiments that demonstrate the performance of the proposed methods. Finally, we conclude with a short summary and some future research topics.

\section{Mathematical representation of randomized neural networks} \label{Sec:RNN}

In this section, we present some notation and mathematical description of randomized neural networks that we use in this paper.

\subsection{Fully connected feedforward neural networks}

A fully connected feedforward neural network consists of applying linear transformations and nonlinear activation functions to the input data in a sequential manner. The nonlinear activation function $\rho(\bx):\Real^d\to\Real^d$ is applied elementwise to a vector $\bx\in\Real^d$. The network has a total of $L+2$ layers, where $L$ is the number of hidden layers. The $l$-th layer has $m_l~(l=0,~\cdots,~L+1)$ neurons. The weights and biases between the $(l-1)$-th and $l$-th layers are denoted by $\bW^{(l)}\in\Real^{m_l\times m_{l-1}}$ and $\bb^{(l)}\in\Real^{m_l\times 1}$, where $l=1,~\cdots,~L+1$.

The fully connected feedforward neural network can be mathematically expressed as follows:
\begin{align} 
    \by^{(0)}&=\bx, \label{Eq:fnn1}\\
     \by^{(l)}&=\rho(\bW^{(l)}\by^{(l-1)}+\bb^{(l)}), \quad l=1,~\cdots,~L, \label{Eq:fnn2}\\
      \bz&=\bW^{(L+1)}\by^{(L)}+\bb^{(L+1)}, \label{Eq:fnn3}
\end{align}
where $\bx\in\Real^{m_0}$ is the input, $\by^{(l)}~(l=1,~\cdots,~L)$ are the hidden layers, and $\bz$ is the output of the network.

\subsection{Randomized neural networks}

With the same structure of a fully connected feedforward neural network, we can generate a randomized neural network (RNN) by randomly assigning the values of some parameters and keeping them fixed during the training process. Specifically, these parameters are the weights and biases of the hidden layers, denoted by $\bcW$ and $\bcB$, respectively, where
\begin{align*}
    \bcW=\set{\bW^{(1)},~\cdots,~\bW^{(L)}}\quad\text{and}\quad\bcB=\set{\bb^{(1)},~\cdots,~\bb^{(L)}}.
\end{align*}
This implies that the output of the last hidden layer, $\by^{(L)}$ in \eqref{Eq:fnn2} is a vector-valued function of the input $\bx$. 
Given a bounded domain $D\subset\Real^{m_0}$, we define $\bphi^{\bcW,~\bcB}_\rho(\bx)=\by^{(L)}$ as a group of basis functions, where each component is a nonlinear function of $\bx$. Then, a randomized neural network function space can be defined as
\begin{align} \label{rnn_space}
    \cN^{\bcW,~\bcB}_{\rho,~m_{L+1}}(D)=\set{\bz(\bx)=\bW\bphi^{\bcW,~\bcB}_\rho(\bx):\bx\in D\subset\Real^{m_0} \quad \forall~\bW\in\Real^{m_{L+1}\times m_L}}.
\end{align}
Note that we set $\bb^{(L+1)}=\bzero$ in the above definition. For convenience, we use $m_0=d$, $m_L=N$ and $m_{L+1}=n$, and we rewrite \eqref{rnn_space} for a fixed $\theta=\{\bcW,~\bcB\}$ as follows:
\begin{align} \label{re_rnn_space}
  [\cN^{\theta}_\rho(D)]^n=\set{\bz(\bx)=\bW\bphi^{\bcW,~\bcB}_\rho(\bx):\bx\in D\subset\Real^{d}\quad \forall~\bW\in\Real^{n\times N}},
\end{align}
which is abbreviated as $\cN^{\theta}_\rho(D)$ when $n=1$.

\section{Local randomized neural networks with hybridized discontinuous Petrov-Galerkin method} \label{Sec:lrnn-hdpg}

In this section, we propose a novel approach that combines local randomized neural networks and a hybridized discontinuous Petrov-Galerkin method to solve the Darcy problem, which models fluid flow in porous media.

We begin with some notation on Sobolev spaces (\cite{Adams1975Sobolev}). Let $k$ be a non-negative integer and $D$ be a bounded domain. The Sobolev space of order $k$ on $D$, denoted by $H^k(D)=W^{k,2}(D)$, consists of functions with square-integrable derivatives up to order $k$. 
The norm and inner product on this space are $\|\cdot\|_{k,D}$ and $(\cdot,~\cdot)_{k,D}$, respectively. In particular, when $D=\Omega$, we simplify the notation above as $\|\cdot\|_k$ and $(\cdot,~\cdot)_k$. If $k=0$, we drop the subscript $k$ for the standard $L^2$ norm and inner product. 

Let $\Omega\subset \Real^d$ be a bounded domain with Lipschitz continuous boundary $\Gamma = \partial\Omega$. We consider the Darcy problem with a Dirichlet boundary condition
\begin{align}
    \bu&=-\underline{\bK}\nabla p\quad&&\text{in }\Omega, \label{Eq:strong-mixed1}\\
    \nabla\cdot\bu&=f\quad&&\text{in }\Omega, \label{Eq:strong-mixed2}\\
    p&=g\quad&&\text{on }\partial\Omega, \label{bd:strong-mixed}
\end{align}
where $p:\Real^d\to\Real$ and $\bu:\Real^d\to\Real^d$ are the unknown pressure and velocity of the flow, respectively. The boundary condition $g\in L^2(\partial\Omega)$ is given. $\underline{\bK}$ is the ratio of the permeability to the viscosity, satisfying $\underline{\bK}=\underline{\bK}^t$. Equation \eqref{Eq:strong-mixed1} implies that the fluid velocity is proportional to the pressure gradient.

Multiplying \eqref{Eq:strong-mixed1} and \eqref{Eq:strong-mixed2} by test functions $\bv$ and $q$, and integrating by parts, we obtain the mixed variational formulation: Find  $(\bu,p)\in H(\text{div},\Omega)\times L^2(\Omega)$ such that
\begin{align}
    (\underline{\bK}^{-1}\bu,~\bv)-(p,~\nabla\cdot\bv)&=-(g,~\bv\cdot\bn)_{\partial\Omega}\quad&&\forall~\bv\in H(\text{div},\Omega),  \label{Eq:weak-mixed1}\\
    (\nabla\cdot\bu,~q)&=(f,~q)\quad&&\forall~q\in L^2(\Omega), \label{Eq:weak-mixed2}
\end{align}
where $H(\text{div},\Omega)=\set{\bv\in[L^2(\Omega)]^d:~\nabla\cdot\bv\in L^2(\Omega)}$ is the space of square-integrable vector fields with square-integrable divergence. Here, $\bn$ is the unit outward normal vector on $\Gamma$.

\subsection{Hybridized discontinuous Galerkin formulation}
\label{subsec:HDG}

We can solve equations \eqref{Eq:weak-mixed1} and \eqref{Eq:weak-mixed2} by using one randomized neural network as in \cite{Shang2023Randomized,Shang2023RNNPG}, however, this approach may be inefficient for some complex problems, such as those with highly oscillatory solutions. Therefore, we adopt the idea of domain decomposition and use one RNN on each subdomain, then apply a hybridized DG scheme to couple neural networks on different subdomains, so that they have better approximation capability.

To introduce the hybridized DG formulation, we present some assumptions and preliminaries. Let $\cT_h$ be a partition of domain $\Omega$ into polygons or polyhedra, such that the elements in $\cT_h$ are closed and simply connected satisfying shape regularity conditions (\cite{Wang2014Weak}). Let $h_K$ denote the diameter of $K\in\cT_h$ and $h=\max_{K\in\cT_h}h_K$. The set of all edges or faces in $\cT_h$ is denoted by $\cE_h$, $\cE_h^b=\set{e\in\cE_h:e\subset\partial\Omega}$ is the set of boundary edges or faces and $\cE_h^i=\cE_h\backslash\cE_h^b$ is the set of interior edges or faces.

We define the following function spaces on $\cT_h$ and $\cE_h$,
\begin{align*}
    \bV&=\set{\bv\in[L^2(\Omega)]^d:~\bv|_K\in[L^2(K)]^d\quad \forall~K\in\cT_h},\\
    \hat{\bV}&=\set{\hat{\bv}\in[L^2(\cE_h)]^d:~\hat{\bv}=\hat{v}\bn_e,~\hat{v}|_e\in L^2(e)\quad\forall~e\in\cE_h},\\
    Q&=\set{q\in L^2(\Omega):~q|_K\in H^1(K)\quad \forall~K\in\cT_h}.
\end{align*}
Here, $\bn_e$ denotes a unit normal vector for $e\in\cE_h$.
Let $(\cdot,~\cdot)_{\cT_h}$, $(\cdot,~\cdot)_{\partial\cT_h}$ and $(\cdot,~\cdot)_{\cE_h^b}$ denote the inner products associated with mesh partition, defined by
\begin{align*}
    (\cdot,~\cdot)_{\cT_h}=\sum_{K\in\cT_h}(\cdot,~\cdot)_K,\quad(\cdot,~\cdot)_{\partial\cT_h}=\sum_{K\in\cT_h}(\cdot,~\cdot)_{\partial K},\quad(\cdot,~\cdot)_{\cE_h^b}=\sum_{e\in\cE_h^b}(\cdot,~\cdot)_e.
\end{align*}

Following the ideas in \cite{Wang2014Weak, Hong2019Unified}, we can derive the hybridized discontinuous Galerkin formulation from \eqref{Eq:weak-mixed1} -- \eqref{Eq:weak-mixed2}: Find $\bu\in\bV$, $\hat{\bu}\in\hat{\bV}$ and $p\in Q$ such that
\begin{align}
    (\underline{\bK}^{-1}\bu,~\bv)_{\cT_h}+(\nabla p,~\bv)_{\cT_h}-(p,~\hat{\bv}\cdot\bn)_{\partial\cT_h}&=-(g,~\hat{\bv}\cdot\bn)_{\cE_h^b}\quad&&\forall~(\bv,\hat{\bv})\in\bV\times\hat{\bV}, \label{Eq:weak-mixedhybrid1}\\
    (\bu,~\nabla q)_{\cT_h}-(\hat{\bu}\cdot\bn,~q)_{\partial\cT_h}&=-(f,~q)_{\cT_h}\quad&&\forall~q\in Q, \label{Eq:weak-mixedhybrid2}
\end{align}
where $(\hat{\bu}\cdot\bn,~q)_{\partial\cT_h} = \sum\limits_{K\in\cT_h} (\hat{\bu}\cdot\bn_K,~q)_{\partial K}$ with $\bn_K$ being the unit outward normal vector on $\partial K$.

\subsection{Discontinuous Petrov-Galerkin method}

To discretize \eqref{Eq:weak-mixedhybrid1} and \eqref{Eq:weak-mixedhybrid2}, we choose suitable function spaces for trial functions $\bu$, $\hat{\bu}$ and $p$, as well as test functions $\bv$, $\hat{\bv}$ and $q$.  
For the trial functions, we use RNN spaces \eqref{re_rnn_space} with different parameters $(\bcW,\bcB,\rho)$. For the test functions, we use piecewise polynomial spaces with different degrees.

Recall that $\Omega\subset\Real^d$ is a bounded domain. We define the trial function spaces as follows:
\begin{align}
    \bV_r&=\set{\bu_r\in[L^2(\Omega)]^d: ~\bu_r|_K\in [\cN_\rho^\theta(K)]^d \quad \forall~K\in\cT_h}, \label{space-u}\\
    \hat{\bV}_r&=\set{\hat{\bu}_r\in [L^2(\cE_h)]^d :~\hat{\bu}_r=\hat{u}_r\bn_e,~\hat{u}_r|_e\in\cN_\rho^\theta(e)\quad\forall~e\in\cE_h}, \label{space-uh}\\
    Q_r&=\set{p_r\in L^2(\Omega): ~p_r|_K\in \cN_\rho^\theta(K)\quad\forall~K\in\cT_h}, \label{space-p}
\end{align}
where $\rho$ and $\theta=\{\bcW,~\bcB\}$ are the same as in \eqref{re_rnn_space}, but they may differ for different edges or faces and elements. Since RNNs are defined on each element $K$ and its edges or faces, we call them local randomized neural network spaces. Usually, we use different weights $\bcW$ and biases $\bcB$ to different elements. However, we can also reduce the computation by using the same $\bcW$ and $\bcB$ for all elements in \eqref{space-u} -- \eqref{space-p}.

\begin{remark}    
The choice of the activation function $\rho$ influences the performance of the method. Different activation functions have different characteristics and advantages. For a complex problem, it may be helpful to use different activation functions for different elements. However, to concentrate on elucidating the idea of our method, we use a uniform activation function for all elements.
\end{remark}

We can choose the following test function spaces:
\begin{align}
    \bV_h&=\set{\bv_h\in[L^2(\Omega)]^d:~\bv_h|_K\in[P_k(K)]^d\quad\forall~K\in\cT_h}, \label{space-v}\\
    \hat{\bV}_h&=\set{\hat{\bv}_h\in [L^2(\cE_h)]^d: ~\hat{\bv}_h=\hat{v}_h\bn_e,~\hat{v}_h|_e\in P_k(e)\quad\forall~e\in \cE_h}, \label{space-vh}\\
    Q_h&=\set{q_h\in L^2(\Omega):~q_h|_K\in P_{k+1}(K)\quad\forall~K\in\cT_h}, \label{space-q}
\end{align}
where $P_k$ denotes the space of polynomials of degree at most $k$.
\begin{remark} \label{Re:deg-poly}
The degree of polynomials in \eqref{space-q} is one degree higher than that in \eqref{space-v} and \eqref{space-vh} to match \eqref{Eq:strong-mixed1}. 
However, the choice of the polynomial degree does not have a significant impact on the results. Therefore, the polynomial degree can be chosen flexibly. 
\end{remark}

\begin{remark} \label{Re:stab-term}
If we use the same function spaces \eqref{space-v} -- \eqref{space-q} for both trial and test functions, the numerical scheme \eqref{Eq:weak-mixedhybrid1} -- \eqref{Eq:weak-mixedhybrid2} may not be stable. In such cases, we usually need to add a stabilization term to ensure the uniqueness of the solution. However, in the proposed method, we solve a least-squares problem that does not require a stabilization term.
\end{remark}

Based on the previous preparation, we present the LRNN-HDPG scheme for \eqref{Eq:strong-mixed1} -- \eqref{bd:strong-mixed}:
Find $\bu_r\in\bV_r$, $\hat{\bu}_r\in\hat{\bV}_r$ and $p_r\in Q_r$ satisfying
\begin{align}
    (\underline{\bK}^{-1}\bu_r,~\bv_h)_{\cT_h}+(\nabla p_r,~\bv_h)_{\cT_h}&=0\quad&&\forall~\bv_h\in\bV_h, \label{AlgoEq:Darcy-eq1}\\
    (p_r,~\hat{\bv}_h\cdot\bn)_{\partial\cT_h}&=(g,~\hat{\bv}_h\cdot\bn)_{\cE_h^b}\quad&&\forall~\hat{\bv}_h\in\hat{\bV}_h, \label{AlgoEq:Darcy-eq2}\\
    (\bu_r,~\nabla q_h)_{\cT_h}-(\hat{\bu}_r\cdot\bn,~q_h)_{\partial\cT_h}&=-(f,~q_h)_{\cT_h}\quad&&\forall~q_h\in Q_h. \label{AlgoEq:Darcy-eq3}
\end{align}

To enhance the stability of the numerical scheme, we can add a stabilization term as in Remark \ref{Re:stab-term}. The stabilized scheme is to find $\bu_r\in\bV_r$, $\hat{\bu}_r\in\hat{\bV}_r$ and $p_r\in Q_r$ that satisfy
\begin{align}
    (\underline{\bK}^{-1}\bu_r,~\bv_h)_{\cT_h}+\eta_h((\bu_r-\hat{\bu}_r)\cdot\bn,~\bv_h\cdot\bn)_{\partial\cT_h}+(\nabla p_r,~\bv_h)_{\cT_h}&=0\quad&&\forall~\bv_h\in\bV_h, \label{AlgoEq:Darcy-stab1}\\
    -\eta_h((\bu_r-\hat{\bu}_r)\cdot\bn,~\hat{\bv}_h\cdot\bn)_{\partial\cT_h}-(p_r,~\hat{\bv}_h\cdot\bn)_{\partial\cT_h}&=-(g,~\hat{\bv}_h\cdot\bn)_{\cE_h^b}\quad&&\forall~\hat{\bv}_h\in\hat{\bV}_h, \label{AlgoEq:Darcy-stab2}\\
    (\bu_r,~\nabla q_h)_{\cT_h}-(\hat{\bu}_r\cdot\bn,~q_h)_{\partial\cT_h}&=-(f,~q_h)_{\cT_h}\quad&&\forall~q_h\in Q_h, \label{AlgoEq:Darcy-stab3}
\end{align}
where $\eta_h(\cdot,~\cdot)_{\partial\cT_h}=\sum\limits_{K\in\cT_h}\eta h_K(\cdot,~\cdot)_{\partial K}$ with a non-negative constant $\eta$. 

Note that \eqref{AlgoEq:Darcy-stab1} -- \eqref{AlgoEq:Darcy-stab3} reduce to \eqref{AlgoEq:Darcy-eq1} -- \eqref{AlgoEq:Darcy-eq3} when $\eta=0$. In the following sections, we use this form with the stabilization term to present the corresponding numerical scheme. As shown in the numerical experiments, the stabilization term is not necessary for the LRNN-HDPG methods. 

\begin{remark} \label{Re:normal-componment}
    Since only $\hat{\bu}_r\cdot\bn$ and $\hat{\bv}_h\cdot\bn$ appear in \eqref{AlgoEq:Darcy-eq1} -- \eqref{AlgoEq:Darcy-eq3} and \eqref{AlgoEq:Darcy-stab1} -- \eqref{AlgoEq:Darcy-stab3}, the trial function space $\hat{\bV}_r$ and the test function space $\hat{\bV}_h$ only involve the normal components of $\hat{\bu}_r$ and $\hat{\bv}_h$.
\end{remark}

\begin{remark}
We can also choose the test function spaces in \eqref{AlgoEq:Darcy-stab1} -- \eqref{AlgoEq:Darcy-stab3} as LRNNs that are the same as the trial function spaces. This leads to an HDG scheme named the LRNN-HDG method. The numerical experiments show that this scheme still achieves high accuracy.
\end{remark}

\subsection{Reduction of degrees of freedom by eliminating $\bu_r$} 
\label{ReduceDOF}

The LRNN-HDPG method can reduce degrees of freedom by applying an elementary transformation of the matrix, similar to other hybridization methods. We use the numerical scheme \eqref{AlgoEq:Darcy-eq1} -- \eqref{AlgoEq:Darcy-eq3} as an example. 

Since $\bu_r$ is locally defined on each element, we can define an operator $\cR$ by \eqref{AlgoEq:Darcy-eq1}. This allows us to replace $\bu_r$ in \eqref{AlgoEq:Darcy-eq3} by $\bu_r=\cR(p_r)$. Then we can get the solutions to  \eqref{AlgoEq:Darcy-eq1} -- \eqref{AlgoEq:Darcy-eq3} by solving a smaller problem: 
Find $\hat{\bu}_r\in\hat{\bV}_r$ and $p_r\in Q_r$ that satisfy
\begin{align}
    (p_r,~\hat{\bv}_h\cdot\bn)_{\partial\cT_h}&=(g,~\hat{\bv}_h\cdot\bn)_{\cE_h^b}\quad&&\forall~\hat{\bv}_h\in\hat{\bV}_h, \label{AlgoEq:Darcy-reduce-eq2}\\
    (\cR(p_r),~\nabla q_h)_{\cT_h}-(\hat{\bu}_r\cdot\bn,~q_h)_{\partial\cT_h}&=-(f,~q_h)_{\cT_h}\quad&&\forall~q_h\in Q_h. \label{AlgoEq:Darcy-reduce-eq3}
\end{align}
And we can get $\bu_r$ by a post-processing step $\bu_r=\cR(p_r)$.

\subsection{A global trace neural network}

The variable $\hat{\bu}$ in \eqref{Eq:weak-mixedhybrid1} -- \eqref{Eq:weak-mixedhybrid2} can be interpreted as a representation of $\bu$ on $\cE_h$, and is used to connect different $\bu$ on two adjacent elements. This means that $\hat{\bu}_r$ in \eqref{AlgoEq:Darcy-eq1} -- \eqref{AlgoEq:Darcy-eq3} couples the LRNNs $\bu_r$ defined on different elements, where $\hat{\bu}_r$ is piecewise LRNNs. In this case, $\hat{\bu}_r$ is defined piecewise on edges or faces, that is, one LRNN is used to approximate $\bu$ on each edge or face. 

In light of the flexibility of the domain of neural network, we can also use a single neural network $\hat{\bu}_\text{global}$ to approximate $\hat{\bu}_r$, that is, this global trace neural network $\hat{\bu}_\text{global}$ is defined on $\overline{\Omega}$, but we only care about its value on edges or faces.

For this purpose, we use $\hat{\bV}_\text{global}$ to replace $\hat{\bV}_r$ in \eqref{space-uh} where $\hat{\bV}_\text{global}$ is defined by
\begin{align} \label{space-uh-global}
    \hat{\bV}_\text{global}=\set{\hat{\bu}_\text{global}:~\hat{\bu}_\text{global}\in [\cN_\rho^\theta(\Omega)]^d}.
\end{align}
The main difference is that they have different dimensions for their domains between $\hat{\bV}_r$ and $\hat{\bV}_\text{global}$, the former is $d-1$ and the latter is $d$. 

Choosing the trial function spaces \eqref{space-u}, \eqref{space-uh-global} and \eqref{space-p}, and the test function spaces \eqref{space-v}, \eqref{space-vh} and \eqref{space-q}, we can also obtain the schemes mentioned above. The numerical experiments show that this scheme works well.

\subsection{Another hybridized DG method}

In Section \ref{subsec:HDG}, we present the hybridized DG method by introducing $\hat{\bu}$ based on \eqref{Eq:weak-mixed1} -- \eqref{Eq:weak-mixed2}. Alternatively, we can propose a different hybridized DG scheme by incorporating $\hat{p}$ (\cite{Cockburn2009Unified}). The corresponding mixed variational formulation seeks $(\bu,p)\in [L^2(\Omega)]^d\times H_g^1(\Omega)$ that satisfies:
\begin{align}
    (\underline{\bK}^{-1}\bu,~\bv)+(\nabla p,~\bv)&=0\quad&&\forall~\bv\in [L^2(\Omega)]^d, \label{Eq:weakII-mixed1}\\
    -(\bu,~\nabla q)&=(f,~q)\quad&&\forall~q\in H_0^1(\Omega), \label{Eq:weakII-mixed2}
\end{align}
where $H_g^1(\Omega)=\set{q:q\in H^1(\Omega),~q|_{\partial\Omega}=g}$. 

To derive the hybridized DG formulation, we define the following function spaces:
\begin{align*}
    \bV_1&=\set{\bv\in[L^2(\Omega)]^d:~\bv|_K\in H(\text{div},K)\quad \forall~K\in\cT_h},\\
    Q_1&=\set{q\in L^2(\Omega):~q|_K\in L^2(K)\quad \forall~K\in\cT_h},\\
    \hat{Q}&=\set{\hat{q}\in L^2(\cE_h):~\hat{q}|_e\in L^2(e)\quad\forall~e\in\cE_h}.
\end{align*}
To address the Dirichlet boundary condition, we introduce two subspaces of $\hat{Q}$,
\begin{align*}
    \hat{Q}^0=\set{\hat{q}\in\hat{Q}:~\hat{q}|_{\partial\Omega}=0}\quad\text{and}\quad\hat{Q}^g=\set{\hat{q}\in\hat{Q}:~\hat{q}|_{\partial\Omega}=g}.
\end{align*}
The discontinuous variational formulation is: Find $\bu\in\bV_1$, $p \in Q_1$ and $\hat{p}\in\hat{Q}^g$ such that
\begin{align}
    (\underline{\bK}^{-1}\bu,~\bv)_{\cT_h}-(p,~\nabla\cdot\bv)_{\cT_h}+(\hat{p},~\bv\cdot\bn)_{\partial\cT_h}&=0\quad&&\forall~\bv\in\bV_1, \label{Eq:weakII-mixedhybrid1}\\
    -(\nabla\cdot\bu,~q)_{\cT_h}+(\bu\cdot\bn,~\hat{q})_{\partial\cT_h}&=-(f,~q)_{\cT_h}\quad&&\forall~(q,\hat{q})\in Q_1\times\hat{Q}^0. \label{Eq:weakII-mixedhybrid2} 
\end{align}

The trial function space $\hat{Q}_r$ for $\hat{p}$ and test function space $\hat{Q}_h$ for $\hat{q}$ can be chosen as follows:
\begin{align}
    \hat{Q}_r^g&=\set{\hat{q}_r\in L^2(\cE_h):~\hat{q}_r|_e\in\cN_\rho^\theta(e),~\hat{q}_r|_{\partial\Omega\cap e}=g\quad\forall~e\in\cE_h}, \label{space-ph}\\
    \hat{Q}_h^0&=\set{\hat{q}_h\in L^2(\cE_h):~\hat{q}_h|_e\in P_{k+1}(e),~\hat{q}_h|_{\partial\Omega\cap e}=0\quad\forall~e\in\cE_h}. \label{space-qh}
\end{align}
Similar to \eqref{AlgoEq:Darcy-stab1} -- \eqref{AlgoEq:Darcy-stab3}, the second LRNN-HDPG scheme for \eqref{Eq:strong-mixed1} -- \eqref{bd:strong-mixed} with a stabilization term is:
Find $\bu_r\in\bV_r$, $p_r\in Q_r$ and $\hat{p}_r\in\hat{Q}_r^g$ satisfying 
\begin{align}
    (\underline{\bK}^{-1}\bu_r,~\bv_h)_{\cT_h}-(p_r,~\nabla\cdot\bv_h)_{\cT_h}+(\hat{p}_r,~\bv_h\cdot\bn)_{\partial\cT_h}&=0\quad&&\forall~\bv_h\in\bV_h,  \label{AlgoEqII:Darcy-stab1}\\
    -(\nabla\cdot\bu_r,~q_h)_{\cT_h}-\tau_h(p_r,~q_h)_{\partial\cT_h}+\tau_h(\hat{p}_r,~q_h)_{\partial\cT_h}&=(f,~q_h)_{\cT_h}\quad&&\forall~q_h\in Q_h,  \label{AlgoEqII:Darcy-stab2}\\
    (\bu_r\cdot\bn,~\hat{q}_h)_{\partial\cT_h}+\tau_h(p_r,~\hat{q}_h)_{\partial\cT_h}-\tau_h(\hat{p}_r,~\hat{q}_h)_{\partial\cT_h}&=0\quad&&\forall~\hat{q}_h\in\hat{Q}_h^0, \label{AlgoEqII:Darcy-stab3}
\end{align}
where $\tau_h(\cdot,~\cdot)_{\partial\cT_h}=\sum\limits_{K\in\cT_h}\tau h_K^{-1}(\cdot,~\cdot)_{\partial K}$ with a non-negative constant $\tau$.

Similarly, when $\tau=0$, an analogous scheme to equations \eqref{AlgoEq:Darcy-eq1} -- \eqref{AlgoEq:Darcy-eq3} is obtained, which remains well-defined by solving a least-squares problem.

\begin{remark}
    To deal with the Dirichlet boundary condition in the space $\hat{Q}_r^g$,
   we introduce a local projection operator $\dQ_h$ that maps $L^2(e)$ onto $\cN_\rho^\theta(e)$ for all $e\in\cE_h^b$. The condition $\hat{p}_r=\dQ_h g$ can be easily imposed on equations \eqref{AlgoEqII:Darcy-stab1} -- \eqref{AlgoEqII:Darcy-stab3}.
\end{remark}

\section{LRNN-HDPG methods for Stokes-Darcy flows}
\label{Sec:Stokes-Darcy-Brinkman}

In this section, we focus on the development of LRNN-HDPG methods tailored for addressing Stokes-Darcy flows and Brinkman equations. We begin by examining the Stokes flow.

\subsection{Stokes equations} \label{Sec:Stokes}

In this subsection, we consider the steady Stokes equations with Dirichlet boundary condition. The approach aligns with the derivations presented in Section \ref{Sec:lrnn-hdpg}. The Stokes equations are
\begin{align}
    -\nabla\cdot(2\nu\underline{\bvarepsilon}(\bu))+\nabla p&=\fb\quad&&\text{in }\Omega,\label{Eq:Stokes-strong-primal}\\
    \nabla\cdot\bu&=0\quad&&\text{in }\Omega,\label{Eq:Stokes-primal-divfree}\\
    \bu&=\bg\quad&&\text{on }\partial\Omega,\label{bd:Stokes-strong-primal}
\end{align}
where $\nu$ represents the fluid viscosity coefficient, $\bu$ denotes the fluid velocity, $p$ is the fluid pressure, and $\underline{\bvarepsilon}(\bu)=(\nabla\bu+(\nabla\bu)^t)/2$ signifies the strain rate. The volume force is given by $\fb$, and the boundary condition $\bg$ satisfies the constraint $\int_{\partial\Omega}\bg\cdot\bn\ds=0$. To ensure a unique solution $(\bu,p)$, the condition $\int_\Omega p\,\dx=0$ is imposed. 

Introducing the stress tensor $\underline{\bsigma}=2\nu\underline{\bvarepsilon}(\bu)-p\underline{\bI}$ where $\underline{\bI}$ is the identity tensor, allows us to reformulate the equations \eqref{Eq:Stokes-strong-primal} -- \eqref{bd:Stokes-strong-primal} as
\begin{align}
    \underline{\bsigma}+p\underline{\bI}-2\nu\underline{\bvarepsilon}(\bu)&=0\quad&&\text{in }\Omega,\label{Eq:Stokes-strong-mixed1}\\
    -\nabla\cdot\underline{\bsigma}&=\fb\quad&&\text{in }\Omega,\label{Eq:Stokes-strong-mixed2}\\
    \nabla\cdot\bu&=0\quad&&\text{in }\Omega,\label{Eq:Stokes-mixed-divfree}\\
    \bu&=\bg\quad&&\text{on }\partial\Omega,\label{bd:Stokes-strong-mixed}\\
    \int_\Omega p\,\dx&=0,\label{Eq:Stokes-strong-intmean}
\end{align}
where $(\underline{\bsigma},\bu,p)$ is the unknown solution. 

Employing the incompressibility condition \eqref{Eq:Stokes-mixed-divfree}, we deduce
\begin{align}
    p=-\frac{1}{d}\text{tr}(\underline{\bsigma}), \label{Eq:Stokes-p-trace}
\end{align}
where $d$ is dimension of $\Omega$. Following \cite{Qian2020Mixed, Qian2023Mixed}, we derive the velocity-stress formulation
\begin{align}
    \underline{\bsigma}^d-2\nu\underline{\bvarepsilon}(\bu)&=0\quad&&\text{in }\Omega,\label{Eq:Stokes-strong-mixed-vp}\\
    -\nabla\cdot\underline{\bsigma}&=\fb\quad&&\text{in }\Omega,\label{Eq:Stokes-mixed-divfree-vp}\\
    \bu&=\bg\quad&&\text{on }\partial\Omega,\label{bd:Stokes-strong-mixed-vp}\\
    \int_\Omega\text{tr}(\underline{\bsigma})\dx&=0,\label{Eq:Stokes-strong-intmean-vp}
\end{align}
where $\underline{\bsigma}^d=\underline{\bsigma}-\frac{1}{d}\text{tr}(\underline{\bsigma})\underline{\bI}$. The primary unknown functions are $\underline{\bsigma}^d$ and $\bu$, with $p$ being computed through post-processing \eqref{Eq:Stokes-p-trace}.

To derive the mixed variational formulation from equations \eqref{Eq:Stokes-strong-mixed-vp} -- \eqref{Eq:Stokes-strong-intmean-vp}, we introduce a function space
\begin{align*}
    [H_0(\text{div},\Omega;\dS)]^d=\set{\underline{\btau}\in[L^2(\Omega)]^{d\times d}:~\underline{\btau}=\underline{\btau}^t,~\nabla\cdot\underline{\btau}\in[L^2(\Omega)]^d,~\int_\Omega\text{tr}(\underline{\btau})\dx=0}.
\end{align*}
By testing equations \eqref{Eq:Stokes-strong-mixed-vp} and \eqref{Eq:Stokes-mixed-divfree-vp} with $\underline{\btau}$ and $\bv$, and applying integration by parts, we arrive at the mixed variational formulation: Find $(\underline{\bsigma},\bu)\in[H_0(\text{div},\Omega;\dS)]^d \times [L^2(\Omega)]^d$ such that
\begin{align}
    \frac{1}{2\nu}(\underline{\bsigma}^d,~\underline{\btau}^d)+(\bu,~\nabla\cdot\underline{\btau})&=(\bg,~\underline{\btau}\bn)_{\partial\Omega}\quad&&\forall~\underline{\btau}\in[H_0(\text{div},\Omega;\dS)]^d,  \label{Eq:Stokes-weak-mixed1}\\
    (\nabla\cdot\underline{\bsigma},~\bv)&=-(\fb,~\bv)\quad&&\forall~\bv\in[L^2(\Omega)]^d, \label{Eq:Stokes-weak-mixed2}
\end{align}

For $\cT_h$ and $\cE_h$, the function spaces are defined as
\begin{align*}
    \underline{\bSigma}&=\set{\underline{\btau}\in[L^2(\Omega)]^{d\times d}:~\underline{\btau}=\underline{\btau}^t,~\underline{\btau}|_K\in[L^2(K)]^{d\times d}\quad\forall~K\in\cT_h,~\int_\Omega\text{tr}(\underline{\btau})\dx=0},\\
    \hat{\underline{\bSigma}}&=\set{\hat{\underline{\btau}}\in[L^2(\cE_h)]^{d\times d}:~\hat{\underline{\btau}}=\hat{\underline{\btau}}^t,~\hat{\underline{\btau}}\bn_e|_e\in[L^2(e)]^d\quad\forall~e\in\cE_h},\\ 
    \bV&=\set{\bv\in[L^2(\Omega)]^d:~\bv|_K\in[H^1(K)]^d\quad\forall~K\in\cT_h}.
\end{align*}
The hybridized discontinuous Galerkin formulation for \eqref{Eq:Stokes-weak-mixed1} -- \eqref{Eq:Stokes-weak-mixed2} seeks $\underline{\bsigma}\in\underline{\bSigma}$, $\hat{\underline{\bsigma}}\in\hat{\underline{\bSigma}}$ and $\bu\in \bV$ such that
\begin{align}
    \frac{1}{2\nu}(\underline{\bsigma}^d,~\underline{\btau}^d)_{\cT_h}-(\underline{\bvarepsilon}(\bu),~\underline{\btau})_{\cT_h}+(\bu,~\hat{\underline{\btau}}\bn)_{\partial\cT_h}&=(\bg,~\hat{\underline{\btau}}\bn)_{\cE_h^b}\quad&&\forall~(\underline{\btau},\hat{\underline{\btau}})\in\underline{\bSigma}\times\hat{\underline{\bSigma}}, \label{Eq:Stokes-weak-mixedhybrid1}\\
    -(\underline{\bsigma},~\underline{\bvarepsilon}(\bv))_{\cT_h}+(\hat{\underline{\bsigma}}\bn,~\bv)_{\partial\cT_h}&=-(\fb,~\bv)_{\cT_h}\quad&&\forall~\bv\in\bV. \label{Eq:Stokes-weak-mixedhybrid2}
\end{align}

As in Section \ref{Sec:lrnn-hdpg}, we discretize the above hybridized DG formulation \eqref{Eq:Stokes-weak-mixedhybrid1} -- \eqref{Eq:Stokes-weak-mixedhybrid2}. Prior to this, we define the trial and test function spaces to be utilized. The trial function spaces are designated as local randomized neural network spaces
\begin{align}
    \underline{\bSigma}_r&=\set{\underline{\bsigma}_r\in[L^2(\Omega)]^{d\times d}:~\underline{\bsigma}_r=\underline{\bsigma}_r^t,~\underline{\bsigma}_r|_K\in [\cN_\rho^\theta(K)]^{d\times d}\quad\forall~K\in\cT_h,~\int_\Omega\text{tr}(\underline{\bsigma_r})\dx=0}, \label{space-sigma}\\
    \hat{\underline{\bSigma}}_r&=\set{\hat{\underline{\bsigma}}_r\in[L^2(\cE_h)]^{d\times d}:~\hat{\underline{\bsigma}}_r=\hat{\underline{\bsigma}}_r^t,~\hat{\underline{\bsigma}}_r\bn_e|_e=\hat{\bsigma}_r\in [\cN_\rho^\theta(e)]^d\quad\forall~e\in\cE_h}, \label{space-sigmah}\\
    \bV_r&=\set{\bu_r\in [L^2(\Omega)]^d:~\bu_r|_K\in [\cN_\rho^\theta(K)]^d \quad\forall~K\in\cT_h}, \label{space-U}
\end{align}
while the test function spaces are conventional piecewise polynomial spaces
\begin{align}
    \underline{\bSigma}_h&=\set{\underline{\btau}_h\in[L^2(\Omega)]^{d\times d}:~\underline{\btau}_h=\underline{\btau}_h^t,~\underline{\btau}_h|_K\in[P_k(K)]^{d\times d}\quad\forall~K\in\cT_h,~\int_\Omega\text{tr}(\underline{\btau_r})\dx=0}, \label{space-tau}\\
    \hat{\underline{\bSigma}}_h&=\set{\hat{\underline{\btau}}_h\in[L^2(\cE_h)]^{d\times d}:~\hat{\underline{\btau}}_h=\hat{\underline{\btau}}_h^t,~\hat{\underline{\btau}}_h\bn_e|_e=\hat{\btau}_h\in[P_k(e)]^d\quad\forall~e\in\cE_h}, \label{space-tauh}\\
    \bV_h&=\set{\bv_h\in[L^2(\Omega)]^d:~\bv_h|_K\in[P_{k+1}(K)]^d\quad\forall~K\in\cT_h}. \label{space-V}
\end{align}
As mentioned in Remark \ref{Re:normal-componment}, the function spaces $\hat{\underline{\bSigma}}_r$ and $\hat{\underline{\bSigma}}_h$ 
share structural similarities with the spaces $\hat{\bV}_r$ and $\hat{\bV}_h$. Once the components of $\hat{\bsigma}_r$ and $\hat{\btau}_h$ are determined, we can choose the components of $\hat{\underline{\bsigma}}_r$ and $\hat{\underline{\btau}}_h$ to fulfill the symmetry conditions of $\hat{\underline{\bsigma}}_r=\hat{\underline{\bsigma}}_r^t$ and $\hat{\underline{\btau}}_h=\hat{\underline{\btau}}_h^t$.

\begin{remark}
    The number of neurons of the output layer in \eqref{space-sigma} is $d(d+1)/2$ because $\underline{\bSigma}_r$ represents the space of $d\times d$ symmetric tensor-valued functions. A symmetric tensor of this dimensionality possesses $d(d+1)/2$ degrees of freedom, which can be represented as a vector of the same length. Consequently, this vector can be approximated by a neural network comprising $d(d+1)/2$ output neurons.
\end{remark}

We present the LRNN-HDPG scheme for equations \eqref{Eq:Stokes-strong-mixed-vp} -- \eqref{Eq:Stokes-strong-intmean-vp} as follows:
Find $\underline{\bsigma}_r\in\underline{\bSigma}_r$, $\hat{\underline{\bsigma}}_r\in\hat{\underline{\bSigma}}_r$ and $\bu_r\in\bV_r$ such that
\begin{align}
    \frac{1}{2\nu}(\underline{\bsigma}_r^d,~\underline{\btau}_h^d)_{\cT_h}+\eta_h((\underline{\bsigma}_r-\hat{\underline{\bsigma}}_r)\bn,~\underline{\btau}_h\bn)_{\partial\cT_h}-(\underline{\bvarepsilon}(\bu_r),~\underline{\btau}_h)_{\cT_h}&=0\quad&&\forall~\underline{\btau}_h\in\underline{\bSigma}_h, \label{AlgoEq:Stokes-stab1}\\
    -\eta_h((\underline{\bsigma}_r-\hat{\underline{\bsigma}}_r)\bn,~\hat{\underline{\btau}}_h\bn)_{\partial\cT_h}+(\bu_r,~\hat{\underline{\btau}}_h\bn)_{\partial\cT_h}&=(\bg,~\hat{\underline{\btau}}_h\bn)_{\cE_h^b}\quad&&\forall~\hat{\underline{\btau}}_h\in\hat{\underline{\bSigma}}_h, \label{AlgoEq:Stokes-stab2}\\
    -(\underline{\bsigma}_r,~\underline{\bvarepsilon}(\bv_h))_{\cT_h}+(\hat{\underline{\bsigma}}_r\bn,~\bv_h)_{\partial\cT_h}&=-(\fb,~\bv_h)_{\cT_h}\quad&&\forall~\bv_h\in\bV_h, \label{AlgoEq:Stokes-stab3}
\end{align}
where $\eta_h(\cdot,~\cdot)_{\partial\cT_h}=\sum\limits_{K\in\cT_h}\eta h_K(\cdot,~\cdot)_{\partial K}$ with $\eta$ being a non-negative constant. 

Note that we add a stabilization term in \eqref{AlgoEq:Stokes-stab1} -- \eqref{AlgoEq:Stokes-stab3}, which is still valid when $\eta=0$ by solving the resulting linear system through least-squares approach.

\subsection{Stokes-Darcy equations} \label{Sec:Stokes-Darcy}

\begin{figure} [!htbp]
    \centering
    \includegraphics[width=0.6\textwidth]{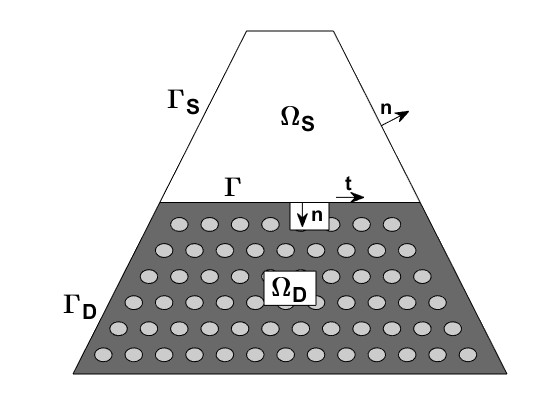}
    \caption{Stokes-Darcy coupling in a 2D domain schematic}
    \label{Fig:Domain-Stokes-Darcy}
\end{figure}

We consider the coupling of Stokes and Darcy flows by an interface $\Gamma$ as shown in Fig. \ref{Fig:Domain-Stokes-Darcy}.
The model consists of two sets of equations that are coupled through interface conditions.

In $\Omega_S$, the flow is governed by the Stokes equations \eqref{Eq:Stokes-strong-primal} -- \eqref{bd:Stokes-strong-primal}, which represent the free flow. We also employ the velocity-stress formulation \eqref{Eq:Stokes-strong-mixed-vp} -- \eqref{Eq:Stokes-strong-intmean-vp} in this subsection. In $\Omega_D$, the porous medium's flow is characterized by the Darcy equations \eqref{Eq:strong-mixed1} -- \eqref{bd:strong-mixed}. The boundaries of $\Omega_S$ and $\Omega_D$ are represented as $\partial\Omega_S=\Gamma_S\cup\Gamma$ and $\partial\Omega_D=\Gamma_D\cup\Gamma$, respectively. We use superscripts $S$ and $D$ to distinguish variables that have the same meaning in different equations, such as $\bu^S$ and $\bu^D$.

The interface conditions are given by
\begin{align}
    \bu^S\cdot\bn&=\bu^D\cdot\bn\quad&&\text{on }\Gamma, \label{Eq:interface1}\\
    (\underline{\bsigma}^S\bn)\cdot\bn&=-p^D\quad&&\text{on }\Gamma, \label{Eq:interface2}\\
    (\underline{\bsigma}^S\bn)\cdot\bt&=-\nu\kappa^{-1}(\bu^S\cdot\bt)\quad&&\text{on }\Gamma, \label{Eq:interface3}
\end{align}
where vectors $\bn$ and $\bt$ are the unit normal and tangent vectors as illustrated in Fig. \ref{Fig:Domain-Stokes-Darcy}.
Here, $\kappa=\frac{\sqrt{(\nu \underline{\bK}\bt)\cdot\bt}}{\alpha}$ is the friction coefficient, with $\underline{\bK}$ and $\nu$ being consistent with those in the Darcy and Stokes equations, respectively. The parameter $\alpha$ is contingent upon the characteristics of $\Gamma$ and is measured by experiment. The first condition \eqref{Eq:interface1} ensures mass conservation, the second condition \eqref{Eq:interface2} equilibrates the normal forces, and the third condition \eqref{Eq:interface3} is the Beavers-Joseph-Saffman law. The second and third conditions can be amalgamated as:
\begin{align}
    \underline{\bsigma}^S\bn+\nu\kappa^{-1}(\bu^S\cdot\bt)\cdot\bt=-p^D\bn\quad\text{on }\Gamma. \label{Eq:interface4}
\end{align}

To enforce the interface conditions \eqref{Eq:interface1} and \eqref{Eq:interface4} for solving Stokes-Darcy equations, many numerical methods use a Lagrange multiplier or an extra constraint. We can adopt similar approaches, however, we use the least-squares method to solve the linear system, which allows us to incorporate additional constraints directly. Thus, we obtain a new linear system by combining \eqref{AlgoEq:Darcy-eq1} -- \eqref{AlgoEq:Darcy-eq3} and \eqref{AlgoEq:Stokes-stab1} -- \eqref{AlgoEq:Stokes-stab3} with additional rows that enforce the interface conditions at select sampling points on $\Gamma$. Let $M$ denote the number of the sampling points. We assume that the partition $\cT_h$ aligns with the interface $\Gamma$. The discrete versions of the interface conditions \eqref{Eq:interface1} -- \eqref{Eq:interface3} are
\begin{align}
     \bu^S_r(\bx_i)\cdot\bn&=\hat{\bu}^D_r(\bx_i)\cdot\bn, \label{Eq:interface-discrete1}\\
     (\hat{\underline{\bsigma}}^S_r(\bx_i)\bn)\cdot\bn&=-p^D_r(\bx_i), \label{Eq:interface-discrete2}\\
     (\hat{\underline{\bsigma}}^S_r(\bx_i)\bn)\cdot\bt&=-\nu\kappa^{-1}(\bx_i)(\bu^S_r(\bx_i)\cdot\bt), \label{Eq:interface-discrete3}
\end{align}
where $\bx_i\in\Gamma$ ($i=1,~\cdots,~M$). These constraints constitute a linear system. The ultimate system is a combination of two diagonal blocks matrices and the linear system \eqref{Eq:interface-discrete1} -- \eqref{Eq:interface-discrete3}, where the two diagonal blocks correspond to \eqref{AlgoEq:Darcy-stab1} -- \eqref{AlgoEq:Darcy-stab3} and \eqref{AlgoEq:Stokes-stab1} -- \eqref{AlgoEq:Stokes-stab3}. 

\begin{remark}
 For simplicity, we choose $\eta=0$ in \eqref{AlgoEq:Darcy-stab1} -- \eqref{AlgoEq:Darcy-stab3} and \eqref{AlgoEq:Stokes-stab1} -- \eqref{AlgoEq:Stokes-stab3}, which means that there is no stabilization term in the proposed method.
\end{remark}

\begin{remark}
    The Beavers-Joseph-Saffman condition \eqref{Eq:interface3} is a simplified version of the following Beavers-Joseph law
    \begin{align}
        (\underline{\bsigma}^S\bn)\cdot\bt&=-\nu\kappa^{-1}((\bu^S-\bu^D)\cdot\bt)\quad{\rm on }~\Gamma. \label{Eq:interface5}
    \end{align}
  The proposed approach can handle both conditions in the same manner. We will verify this for the general Beavers-Joseph condition in the numerical experiments.
\end{remark}

\subsection{Brinkman equations} \label{Sec:Brinkman}

In this subsection, we present another approach of coupling Stokes flow and Darcy flow, leading to the Brinkman equations (\cite{Brinkman1949Calculation}). The Brinkman equations are given by
\begin{align}
    \nu\underline{\bkappa}^{-1}\bu-\nabla\cdot(2\nu\underline{\bvarepsilon}(\bu))+\nabla p&=\fb\quad&&\text{in }\Omega,\label{Eq:Brinkman-strong-primal}\\
    \nabla\cdot\bu&=0\quad&&\text{in }\Omega,\label{Eq:Brinkman-primal-divfree}\\
    \bu&=\bg\quad&&\text{on }\partial\Omega,\label{bd:Brinkman-strong-primal}
\end{align}
where the tensor $ \underline{\bkappa}$ represents the permeability of the fluid, and the other parameters are consistent with those in Section \ref{Sec:Stokes}. The Brinkman equations can be interpreted as the Stokes equations with an additional term $\nu\underline{\bkappa}^{-1}\bu$. The flow is dominated by Darcy's law when $\bkappa$ approaches zero, and by Stokes flow when $\bkappa$ is large.

To solve for $(\bu,p)$, we empoly the velocity-stress formulation for \eqref{Eq:Brinkman-strong-primal} -- \eqref{bd:Brinkman-strong-primal}, which is
\begin{align}
    \underline{\bsigma}^d-2\nu\underline{\bvarepsilon}(\bu)&=0\quad&&\text{in }\Omega,\label{Eq:Brinkman-strong-mixed-vp}\\
    \nu\underline{\bkappa}^{-1}\bu-\nabla\cdot\underline{\bsigma}&=\fb\quad&&\text{in }\Omega,\label{Eq:Brinkman-mixed-divfree-vp}\\
    \bu&=\bg\quad&&\text{on }\partial\Omega,\label{bd:Brinkman-strong-mixed-vp}\\
    \int_\Omega\text{tr}(\underline{\bsigma})\dx&=0.\label{Eq:Brinkman-strong-intmean-vp}
\end{align}
The symbols in this section have the same meanings as in Section \ref{Sec:Stokes}. We omit the derivation and present the numerical scheme for \eqref{Eq:Brinkman-strong-mixed-vp} -- \eqref{Eq:Brinkman-strong-intmean-vp}: 
Find $\underline{\bsigma}_r\in\underline{\bSigma}_r$, $\hat{\underline{\bsigma}}_r\in\hat{\underline{\bSigma}}_r$ and $\bu_r\in\bV_r$ satisfying
\begin{align}
    \frac{1}{2\nu}(\underline{\bsigma}_r^d,~\underline{\btau}_h^d)_{\cT_h}+\eta_h((\underline{\bsigma}_r-\hat{\underline{\bsigma}}_r)\bn,~\underline{\btau}_h\bn)_{\partial\cT_h}-(\underline{\bvarepsilon}(\bu_r),~\underline{\btau}_h)_{\cT_h}&=0\quad&&\forall~\underline{\btau}_h\in\underline{\bSigma}_h, \label{AlgoEq:Brinkman-stab1}\\
    -\eta_h((\underline{\bsigma}_r-\hat{\underline{\bsigma}}_r)\bn,~\hat{\underline{\btau}}_h\bn)_{\partial\cT_h}+(\bu_r,~\hat{\underline{\btau}}_h\bn)_{\partial\cT_h}&=(\bg,~\hat{\underline{\btau}}_h\bn)_{\cE_h^b}\quad&&\forall~\hat{\underline{\btau}}_h\in\hat{\underline{\bSigma}}_h, \label{AlgoEq:Brinkman-stab2}\\
    -\nu(\underline{\bkappa}^{-1}\bu_r,~\bv_h)_{\cT_h}-(\underline{\bsigma}_r,~\underline{\bvarepsilon}(\bv_h))_{\cT_h}+(\hat{\underline{\bsigma}}_r\bn,~\bv_h)_{\partial\cT_h}&=-(\fb,~\bv_h)_{\cT_h}\quad&&\forall~\bv_h\in\bV_h, \label{AlgoEq:Brinkman-stab3}
\end{align}
where $\eta_h(\cdot,~\cdot)_{\partial\cT_h}=\sum\limits_{K\in\cT_h}\eta h_K(\cdot,~\cdot)_{\partial K}$ with a non-negative constant $\eta$.

\section{Numerical experiments} \label{Sec:NumerExper}

In this section, we present various numerical experiments to demonstrate the effectiveness of the proposed methods.

We use RNNs with a single hidden layer for all the numerical examples. We use the hyperbolic tangent function $\text{tanh}(x)=\frac{e^x-e^{-x}}{e^x+e^{-x}}$ as the activation function.
This gives us the basis functions
\begin{align*}
    \bphi_\rho^{\bcW,~\bcB}(\bx)=\text{tanh}(\bW^{(1)}\bx+\bb^{(1)})
\end{align*}
on each element, and
\begin{align*}
    \bphi_\rho^{\bcW,~\bcB}(t)=\text{tanh}(\bW^{(1)}t+\bb^{(1)})+\text{tanh}(\bW^{(1)}_{\text{flip}}(1-t)+\bb^{(1)}_{\text{flip}}),
\end{align*}
on each edge,
where $\bW^{(1)}=(w_1,~\cdots,~w_N)^t$ and $\bW^{(1)}_{\text{flip}}=(w_N,~\cdots,~w_1)^t$. This special structure is a natural consequence of our method and is very convenient for implementing randomized neural networks. Since the weights of the nonlinear part of the randomized neural networks are fixed, we can precisely differentiate the basis functions without using automatic differentiation or finite difference schemes. 

For the least-squares solver, we apply the QR factorization with column pivoting in all the numerical examples. Because of the randomness of $\bW^{(1)}$ and $\bb^{(1)}$, all the numerical results are the averages of multiple experiments. Unless otherwise stated, we perform 10 experiments for each case.

For the reader's convenience, we summarize the main notation used in this paper, some of which have been introduced before. 
\begin{itemize}
    \item Let $N$ be the number of neurons in the hidden layer. We use subscripts to specify the number of neurons for each RNN function. For example, $N_{\bu_r}$, $N_{\hat{\bu}_r}$ and $N_{p_r}$ denote the number of neurons for $\bu_r$, $\hat{\bu}_r$ and $p_r$, respectively.
    \item $d$ is the number of input layer neurons as well as the dimension of the domain. We only consider the case $d=2$, which means the domain $\Omega$ is two-dimensional. 
    \item $k$ is the degree of the polynomial. We have $\text{dim}(P_k)=k+1$ and $\text{dim}([P_k]^2)=(k+1)(k+2)/2$. 
    \item We randomly generate each element of $\bW^{(1)}$ and $\bb^{(1)}$ from the uniform distribution $\cU(-r,~r)$. 
    \item DoF is the degree of freedom as well as the number of columns of the global stiffness matrix. 
    \item We use relative $L^2$ and semi-$H^1$ errors defined by 
    \begin{align*}
        e_0(\ast)=\frac{\|\ast_r-\ast_{\text{real}}\|}{\|\ast_{\text{real}}\|}\quad\text{and}\quad e_1(\ast)=\frac{|\ast_r-\ast_{\text{real}}|_{H^1}}{|\ast_{\text{real}}|_{H^1}}.
    \end{align*}
    In Example \ref{ex2}, we also use relative $L^1$ error, which is defined by 
    \begin{align*}
        \epsilon_1(\ast)=\frac{\|\ast_r-\ast_{\text{real}}\|_{L^1}}{\|\ast_{\text{real}}\|_{L^1}}.
    \end{align*}
    \item For all problems, we use uniform rectangular meshes, which divide the domain into $N_x\times N_y$ equal-sized rectangles. Mesh size is denoted by $h=N_x^{-1}=N_y^{-1}$. 
\end{itemize}

\begin{example}[Darcy flows through anisotropic porous media] \label{ex1}
    Let $\Omega=(0,~1)^2$ and $\underline{\bK}=\left(\begin{matrix}10 & 2 \\ 2 & 100 \end{matrix}\right)$, we consider the Darcy problem \eqref{Eq:strong-mixed1} -- \eqref{bd:strong-mixed}. The exact solution is given by
    \begin{align*}
        p(x,~y)=x(1-x)y(1-y)e^{xy},
    \end{align*}
and we can directly calculate the source term $f$ and set Dirichlet condition $g=p$ on the boundary $\partial\Omega$.
\end{example}

To match the dimensions of the polynomial spaces, we choose different numbers of neurons for $\bu_r$, $\hat{\bu}_r$ and $p_r$, such that $N_{\bu_r}=\text{dim}([P_k]^2)$, $N_{\hat{\bu}_r}=\text{dim}(P_k)$, and $N_{p_r}=\text{dim}([P_{k+1}]^2)$. As stated in Remark \ref{Re:deg-poly}, other choices are also feasible.

\begin{table}[!htpb]
\centering
\begin{tabular}{|c||c|c|c|c||c|c|c|c||}
\hline
& \multicolumn{4}{c||}{$N_{\bu_r}=6$, $N_{\hat{\bu}_r}=3$, $N_{p_r}=10$, DoF $=270$} & \multicolumn{4}{c||}{$N_{\bu_r}=15$, $N_{\hat{\bu}_r}=5$, $N_{p_r}=21$, DoF $=579$} \\ 
\hline
& \multicolumn{3}{c|}{$\eta=0$} & $\eta=1$ & \multicolumn{3}{c|}{$\eta=0$} & $\eta=1$ \\ 
\hline
\diagbox{$k$}{err} & $e_0(p)$ & $e_1(p)$ & $e_0(\bu)$ & $e_0(p)$ & $e_0(p)$ & $e_1(p)$ & $e_0(\bu)$ & $e_0(p)$ \\
\hline                       
3 & 3.29e-02 & 1.32e-01 & 1.74e-01 & 1.11e+00 & 1.67e-02 & 1.09e-01 & 4.20e-02 & 3.34e-01 \\
\hline
4 & 4.67e-02 & 1.82e-01 & 1.07e-01 & 1.06e+00 & 4.89e-04 & 3.69e-03 & 3.41e-03 & 4.94e-01 \\
\hline
5 & 5.57e-02 & 2.02e-01 & 1.30e-01 & 8.93e-01 & 4.27e-04 & 3.09e-03 & 2.89e-03 & 1.81e-02 \\
\hline
6 & 6.57e-02 & 2.28e-01 & 1.25e-01 & 9.05e-01 & 7.24e-04 & 4.91e-03 & 3.45e-03 & 2.08e-02 \\
\hline
7 & 7.47e-02 & 2.45e-01 & 1.44e-01 & 1.09e+00 & 6.58e-04 & 4.32e-03 & 3.18e-03 & 1.64e-02 \\
\hline
8 & 7.32e-02 & 2.40e-01 & 1.29e-01 & 1.07e+00 & 6.52e-04 & 4.22e-03 & 3.15e-03 & 1.76e-02 \\
\hline
\hline
& \multicolumn{4}{c||}{$N_{\bu_r}=28$, $N_{\hat{\bu}_r}=7$, $N_{p_r}=36$, DoF $=996$} & \multicolumn{4}{c||}{$N_{\bu_r}=45$, $N_{\hat{\bu}_r}=9$, $N_{p_r}=55$, DoF $=1521$} \\
\hline
& \multicolumn{3}{c|}{$\eta=0$} & $\eta=1$ & \multicolumn{3}{c|}{$\eta=0$} & $\eta=1$ \\ 
\hline
\diagbox{$k$}{err} & $e_0(p)$ & $e_1(p)$ & $e_0(\bu)$ & $e_0(p)$ & $e_0(p)$ & $e_1(p)$ & $e_0(\bu)$ & $e_0(p)$ \\
\hline    
3 & 1.52e-02 & 1.03e-01 & 4.12e-02 & 3.04e-01 & 1.29e-02 & 8.73e-02 & 3.36e-02 & 3.05e-01 \\
\hline
4 & 3.45e-03 & 2.65e-02 & 8.74e-03 & 5.68e-02 & 4.54e-03 & 3.46e-02 & 1.12e-02 & 5.08e-02 \\
\hline
5 & 2.24e-04 & 2.05e-03 & 1.03e-03 & 9.02e-03 & 3.21e-04 & 2.89e-03 & 1.51e-03 & 8.51e-03 \\
\hline
6 & 3.74e-05 & 3.78e-04 & 2.75e-04 & 8.96e-04 & 5.04e-05 & 5.50e-04 & 2.60e-04 & 1.16e-03 \\
\hline
7 & 4.33e-05 & 4.18e-04 & 2.00e-04 & 8.98e-04 & 3.39e-05 & 3.72e-04 & 1.18e-04 & 8.47e-04 \\
\hline
8 & 4.54e-05 & 4.47e-04 & 2.03e-04 & 1.16e-03 & 3.59e-05 & 3.83e-04 & 1.42e-04 & 4.58e-04 \\
\hline
\end{tabular}
\caption{Relative errors of LRNN-HDPG method \eqref{AlgoEq:Darcy-stab1} -- \eqref{AlgoEq:Darcy-stab3} for various numbers of neurons and polynomial degrees on a fixed mesh of size $h=3^{-1}$ in Example \ref{ex1}.}
\label{table:ex1-neuron-mesh33}
\end{table}

\begin{table}[!htpb]
\centering
\begin{tabular}{|c||c|c|c||c|c|c||}
\hline
& \multicolumn{3}{c||}{$N_{\bu_r}=6$, $N_{\hat{\bu}_r}=3$} & \multicolumn{3}{c||}{$N_{\bu_r}=15$, $N_{\hat{\bu}_r}=5$} \\ 
& \multicolumn{3}{c||}{$N_{p_r}=10$, DoF $=162\ (270)$} & \multicolumn{3}{c||}{$N_{p_r}=21$, DoF $=309\ (579)$} \\ 
\hline
\diagbox{$k$}{err} & $e_0(p)$ & $e_1(p)$ & $e_0(\bu)$ & $e_0(p)$ & $e_1(p)$ & $e_0(\bu)$ \\
\hline                       
3 & 4.56e-02 & 2.06e-01 & 1.18e-01 & 1.14e-02 & 7.42e-02 & 5.20e-02 \\
\hline
4 & 5.57e-02 & 2.11e-01 & 8.51e-02 & 3.05e-04 & 2.36e-03 & 3.69e-03 \\
\hline
5 & 7.50e-02 & 2.95e-01 & 1.08e-01 & 5.26e-04 & 3.50e-03 & 2.07e-03 \\
\hline
6 & 7.00e-02 & 2.59e-01 & 1.15e-01 & 5.38e-04 & 2.74e-03 & 1.30e-03 \\
\hline
7 & 8.98e-02 & 3.17e-01 & 1.21e-01 & 8.37e-04 & 3.45e-03 & 1.68e-03 \\
\hline
8 & 7.40e-02 & 2.63e-01 & 1.15e-01 & 1.10e-03 & 4.07e-03 & 1.90e-03 \\
\hline
\hline
& \multicolumn{3}{c||}{$N_{\bu_r}=28$, $N_{\hat{\bu}_r}=7$} & \multicolumn{3}{c||}{$N_{\bu_r}=45$, $N_{\hat{\bu}_r}=9$} \\
& \multicolumn{3}{c||}{$N_{p_r}=36$, DoF $=492\ (996)$} & \multicolumn{3}{c||}{$N_{p_r}=55$, DoF $=711\ (1521)$} \\
\hline
\diagbox{$k$}{err} & $e_0(p)$ & $e_1(p)$ & $e_0(\bu)$ & $e_0(p)$ & $e_1(p)$ & $e_0(\bu)$ \\
\hline    
3 & 1.02e-02 & 6.96e-02 & 3.57e-02 & 1.21e-02 & 8.58e-02 & 6.54e-02 \\
\hline
4 & 8.35e-04 & 6.49e-03 & 4.62e-03 & 1.48e-03 & 1.16e-02 & 7.56e-03 \\
\hline
5 & 1.05e-04 & 1.01e-03 & 5.22e-04 & 1.33e-04 & 1.28e-03 & 6.87e-04 \\
\hline
6 & 2.06e-05 & 1.81e-04 & 1.77e-04 & 6.00e-05 & 6.17e-04 & 1.40e-04 \\
\hline
7 & 4.77e-05 & 4.45e-04 & 1.19e-04 & 3.39e-05 & 3.37e-04 & 7.93e-05 \\
\hline
8 & 7.91e-05 & 6.21e-04 & 1.59e-04 & 1.10e-05 & 7.17e-05 & 3.58e-05 \\
\hline
\end{tabular}
\caption{Relative errors of the degrees of freedom reduction scheme \eqref{AlgoEq:Darcy-reduce-eq2} -- \eqref{AlgoEq:Darcy-reduce-eq3} for various numbers of neurons and polynomial degrees on a fixed mesh of size $h=3^{-1}$ in Example \ref{ex1}.}
\label{table:ex1-neuron-mesh-reduce33}
\end{table}

\begin{table}[!htpb]
\centering
\begin{tabular}{|c||c|c|c||c|c|c||}
\hline
& \multicolumn{3}{c||}{$N_{\bu_r}=6$, $N_{\hat{\bu}_r}=72$} & \multicolumn{3}{c||}{$N_{\bu_r}=15$, $N_{\hat{\bu}_r}=120$} \\ 
& \multicolumn{3}{c||}{$N_{p_r}=10$, DoF $=270$} & \multicolumn{3}{c||}{$N_{p_r}=21$, DoF $=579$} \\ 
\hline
\diagbox{$k$}{err} & $e_0(p)$ & $e_1(p)$ & $e_0(\bu)$ & $e_0(p)$ & $e_1(p)$ & $e_0(\bu)$ \\
\hline                       
3 & 3.34e-02 & 1.32e-01 & 1.22e-01 & 1.95e-01 & 1.50e+00 & 2.96e-02 \\
\hline
4 & 4.55e-02 & 1.68e-01 & 1.06e-01 & 2.83e-04 & 2.18e-03 & 2.36e-03 \\
\hline
5 & 6.35e-02 & 2.19e-01 & 1.20e-01 & 7.06e-04 & 5.50e-03 & 3.06e-03 \\
\hline
6 & 5.33e-02 & 1.81e-01 & 1.02e-01 & 7.03e-04 & 4.97e-03 & 2.20e-03 \\
\hline
7 & 8.08e-02 & 2.48e-01 & 1.15e-01 & 1.09e-03 & 7.12e-03 & 3.34e-03 \\
\hline
8 & 8.34e-02 & 2.51e-01 & 1.25e-01 & 1.48e-03 & 8.24e-03 & 3.02e-03 \\
\hline
\hline
& \multicolumn{3}{c||}{$N_{\bu_r}=28$, $N_{\hat{\bu}_r}=168$} & \multicolumn{3}{c||}{$N_{\bu_r}=45$, $N_{\hat{\bu}_r}=216$} \\
& \multicolumn{3}{c||}{$N_{p_r}=36$, DoF $=996$} & \multicolumn{3}{c||}{$N_{p_r}=55$, DoF $=1521$} \\
\hline
\diagbox{$k$}{err} & $e_0(p)$ & $e_1(p)$ & $e_0(\bu)$ & $e_0(p)$ & $e_1(p)$ & $e_0(\bu)$ \\
\hline    
3 & 1.89e-01 & 1.46e+00 & 4.39e-02 & 4.89e-01 & 3.78e+00 & 1.18e-01 \\
\hline
4 & 6.06e-02 & 5.28e-01 & 1.69e-02 & 7.71e-02 & 6.68e-01 & 1.52e-02 \\
\hline
5 & 5.70e-04 & 6.20e-03 & 1.40e-03 & 2.20e-03 & 2.48e-02 & 3.39e-03 \\
\hline
6 & 3.30e-04 & 3.42e-03 & 9.57e-04 & 8.25e-04 & 9.73e-03 & 2.30e-03 \\
\hline
7 & 1.86e-04 & 1.79e-03 & 3.95e-04 & 3.21e-04 & 3.45e-03 & 6.61e-04 \\
\hline
8 & 2.64e-04 & 2.25e-03 & 4.78e-04 & 3.25e-04 & 3.39e-03 & 6.02e-04 \\
\hline
\end{tabular}
\caption{Relative errors of LRNN-HDPG method with global trace neural network \eqref{space-uh-global} for various numbers of neurons and polynomial degrees on a fixed mesh of size $h=3^{-1}$ in Example \ref{ex1}.}
\label{table:ex1-neuron-mesh-global33}
\end{table}

We use a $3\times3$ uniform rectangular partition and choose the uniform distribution $\cU(-0.6,~0.6)$. Table \ref{table:ex1-neuron-mesh33} shows relative $L^2$ errors for $\bu_r$ and $p_r$, relative semi-$H^1$ errors for $p_r$ when $\eta=0$, and the relative $L^2$ errors for $p_r$ when $\eta=1$. Here, $\eta$ is the constant in the scheme \eqref{AlgoEq:Darcy-stab1} -- \eqref{AlgoEq:Darcy-stab3}. 
The numerical results indicate that the LRNN-HDPG method achieves better accuracy with an increased number of degrees of freedom and test functions. As demonstrated in Table \ref{table:ex1-neuron-mesh33}, the stabilization term does not enhance the results; hence, we exclusively evaluate the LRNN-HDPG methods with $\eta=0$ in subsequent numerical experiments.

Next, we solve this problem by the degrees of freedom reduction scheme \eqref{AlgoEq:Darcy-reduce-eq2} -- \eqref{AlgoEq:Darcy-reduce-eq3}. All parameters remain consistent with those in Table \ref{table:ex1-neuron-mesh33}, except for the total degrees of freedom. 
Table \ref{table:ex1-neuron-mesh-reduce33} shows that all relative errors remain the precision of those in Table \ref{table:ex1-neuron-mesh33}, validating the effectiveness of this scheme. It's worth noting that the numbers in parentheses after DoF in Table \ref{table:ex1-neuron-mesh-reduce33} represent the degrees of freedom prior to reduction.

Subsequently, we assess the LRNN-HDPG method incorporating the global trace neural network defined in \eqref{space-uh-global}. The number of neurons of RNN in $\hat{\bV}_{\text{global}}$ is denoted by $N_{\hat{\bu}_r}$. The parameters are aligned with those presented in Table \ref{table:ex1-neuron-mesh33}. Since there is only one neural network in $\hat{\bV}_{\text{global}}$, $N_{\hat{\bu}_r}$ is chosen such that the total degrees of freedom for this scheme are same with those in Table \ref{table:ex1-neuron-mesh33}.
According to Table \ref{table:ex1-neuron-mesh-global33}, the LRNN-HDPG method, when applied with the global trace neural network, yields accuracy comparable to the results documented in Table \ref{table:ex1-neuron-mesh33}.

\begin{figure}[!htbp]
    \centering
    \begin{subfigure}[b]{0.32\textwidth}
        \centering
        \includegraphics[width=\textwidth]{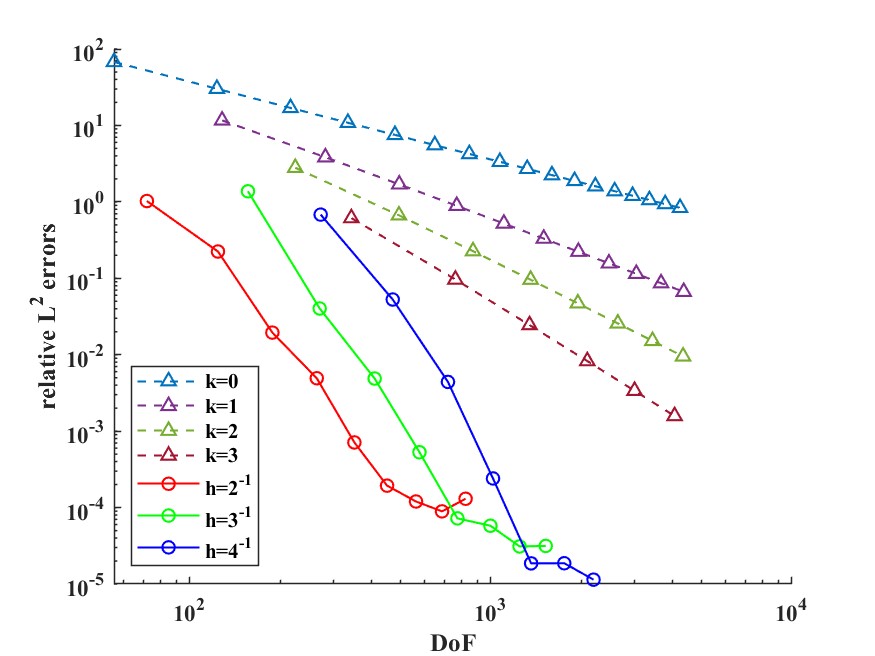}
        \caption{Relative $L^2$ errors of $p_r$}
        \label{fig:ex1-L2p}
    \end{subfigure}
    \begin{subfigure}[b]{0.32\textwidth}
        \centering
        \includegraphics[width=\textwidth]{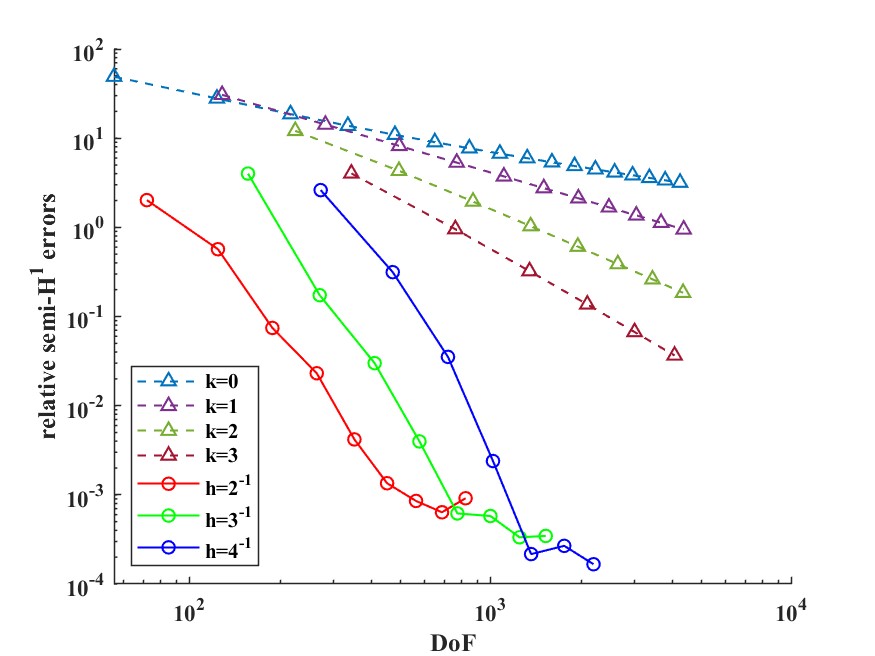}
        \caption{Relative semi-$H^1$ errors of $p_r$}
        \label{fig:ex1-H1p}
    \end{subfigure}
    \begin{subfigure}[b]{0.32\textwidth}
        \centering
        \includegraphics[width=\textwidth]{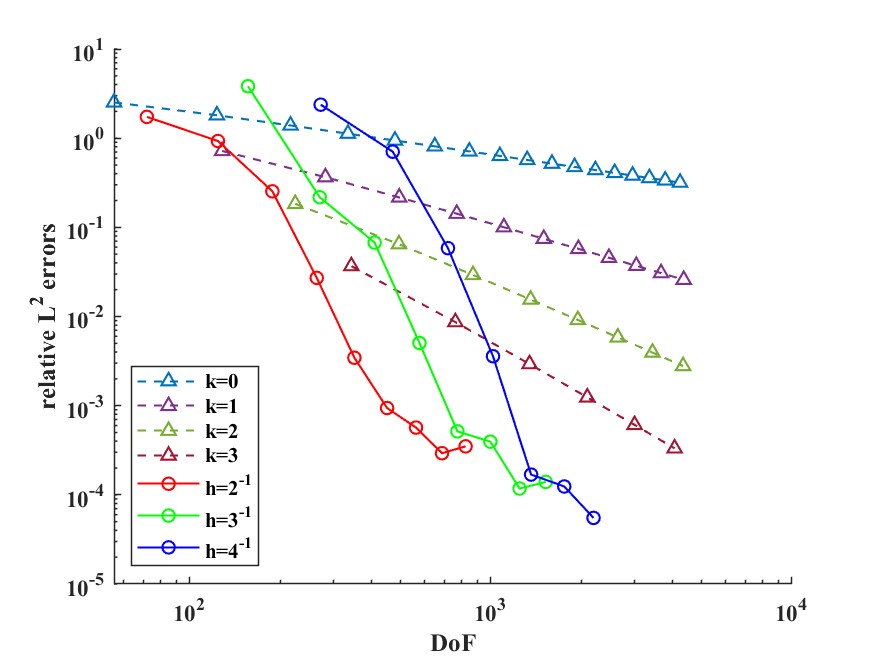}
        \caption{Relative $L^2$ errors of $\bu_r$}
        \label{fig:ex1-L2u}
    \end{subfigure}
    \caption{Comparison of LRNN-HDPG method \eqref{AlgoEq:Darcy-stab1} -- \eqref{AlgoEq:Darcy-stab3} (solid lines) and WG-MFEM (dotted lines) with $([P_k]^2,~P_k,~P_{k+1})$ in Example \ref{ex1}.} 
    \label{fig:ex1}
\end{figure}

Finally, we compare our method with the weak Galerkin mixed finite element method (WG-MFEM) proposed in \cite{Wang2013Weak}. 
With a fixed partition, the degree of freedom is modulated by varying the neuron count in the LRNN-HDPG method, whereas WG-MFEM alters the mesh size for a fixed polynomial degree. Figure \ref{fig:ex1} illustrates the errors across different polynomial degrees and mesh sizes. Corresponding to different mesh sizes $h$, distinct uniform distributions $\cU(-r,~r)$ are applied, where $r=0.8$, $0.6$ and $0.8$ for $h=2^{-1}$, $3^{-1}$ and $4^{-1}$, respectively. We observe that the LRNN-HDPG method can achieve highly accurate numerical solutions with fewer degrees of freedom.

\begin{example}[Darcy problem with highly oscillatory solutions]
 \label{ex2}
    Let $\Omega=(0,~1)^2$, $\underline{\bK}=\underline{\bI}_{2\times2}$, and the exact solution be 
    \begin{align*}
        p(x,~y)=\frac{1}{m}\sum_{i=1}^m\sin(2^i\pi x)\sin(2^i\pi y),
    \end{align*}
    where $m$ is a positive integer. This example is taken from \cite{Dolean2023Multilevel}.
\end{example}

We choose the same RNN spaces for both test and trial function spaces, resulting in a square stiffness matrix. This leads to the LRNN-HDG method. We can also reduce the degrees of freedom by using the approach in Section \ref{ReduceDOF}. We apply the reduction scheme \eqref{AlgoEq:Darcy-reduce-eq2} -- \eqref{AlgoEq:Darcy-reduce-eq3} in this numerical experiment.

\begin{table}[!htbp]
\centering
\begin{tabular}{|c|c|c|c|c|c|c|c|c|}
\hline
& \multicolumn{2}{c|}{$h=2^{-1}$, $r=2.1$} & \multicolumn{2}{c|}{$h=4^{-1}$, $r=1.9$} & \multicolumn{2}{c|}{$h=6^{-1}$, $r=1.6$} & \multicolumn{2}{c|}{$h=8^{-1}$, $r=1.6$} \\ 
\hline
\diagbox{$k_0$}{err} & $e_0(p)$ & $\epsilon_1(p)$ & $e_0(p)$ & $\epsilon_1(p)$ & $e_0(p)$ & $\epsilon_1(p)$ & $e_0(p)$ & $\epsilon_1(p)$ \\
\hline
2 & 9.73e-01 & 9.62e-01 & 8.48e-01 & 8.43e-01 & 5.25e-01 & 5.35e-01 & 3.61e-01 & 3.70e-01 \\ 
\hline
3 & 8.35e-01 & 7.69e-01 & 5.85e-01 & 6.03e-01 & 2.30e-01 & 2.31e-01 & 9.81e-02 & 9.51e-02 \\ 
\hline
4 & 7.44e-01 & 7.24e-01 & 4.39e-01 & 4.51e-01 & 7.19e-02 & 7.16e-02 & 1.76e-02 & 1.71e-02 \\ 
\hline
5 & 5.75e-01 & 6.07e-01 & 1.56e-01 & 1.56e-01 & 1.94e-02 & 1.92e-02 & 4.51e-03 & 4.24e-03 \\ 
\hline
6 & 5.51e-01 & 5.74e-01 & 6.47e-02 & 6.19e-02 & 5.45e-03 & 5.07e-03 & 8.53e-04 & 7.93e-04 \\ 
\hline
7 & 5.40e-01 & 5.58e-01 & 2.39e-02 & 2.27e-02 & 1.74e-03 & 1.61e-03 & 3.17e-04 & 2.72e-04 \\ 
\hline
8 & 4.95e-01 & 5.04e-01 & 1.04e-02 & 9.60e-03 & 8.83e-04 & 7.84e-04 & 1.62e-04 & 1.32e-04 \\ 
\hline
9 & 4.04e-01 & 4.11e-01 & 6.11e-03 & 5.81e-03 & 5.00e-04 & 4.63e-04 & 9.95e-05 & 8.18e-05 \\ 
\hline
10 & 3.21e-01 & 3.21e-01 & 4.67e-03 & 4.26e-03 & 4.07e-04 & 3.90e-04 & 7.87e-05 & 6.12e-05 \\
\hline
\end{tabular}
\caption{Relative $L^2$ and $L^1$ errors of LRNN-HDG method for various numbers of neurons and mesh sizes when $m=3$ in Example \ref{ex2}. }
\label{table:ex2-HDG-n3}
\end{table}

\begin{table}[!htbp]
\centering
\begin{tabular}{|c|c|c|c|c|c|c|c|c}
\hline
$(m,~h^{-1},~r)$ & $(1,~8,~0.9)$ & $(2,~8,~1.1)$ & $(3,~16,~1.1)$ & $(4,~16,~1.2)$ & $(5,~32,~1.1)$ & $(6,~32,~1.3)$ \\
\hline
$e_0(p)$ & 7.1819e-06 & 9.6998e-05 & 1.0392e-04 & 3.9073e-03 & 3.8311e-03 & 1.0291e-01 \\
\hline
$\epsilon_1(p)$ & 6.1324e-06 & 8.7189e-05 & 9.0233e-05 & 3.5010e-03 & 3.4167e-03& 9.7048e-02 \\
\hline
\end{tabular}
\caption{Relative $L^2$ and $L^1$ errors of LRNN-HDG method with $k_0=5$ for different $m$ in Example \ref{ex2}. }
\label{table:ex2-HDG-n123456}
\end{table}

We choose several uniform rectangle partitions to solve the problem for $m=3$. The number of neurons satisfies 
\begin{align*}
    N_{\bu_r}=\frac{1}{2}(k_0+1)(k_0+2),\quad N_{\hat{\bu}_r}=k_0+1,\quad N_{p_r}=\frac{1}{2}(k_0+2)(k_0+3).
\end{align*}
For different mesh sizes $h$, we have different uniform distributions $\cU(-r,~r)$. Relative $L^2$ and $L^1$ errors are shown in Table \ref{table:ex2-HDG-n3} for different $k_0$ and $h$. We observe that the LRNN-HDG method improves the accuracy with more degrees of freedom and smaller $h$.

We repeat this example for $m=1,\cdots,6$. We can choose a fine enough uniform mesh size to ensure accuracy. However, a coarse mesh is sufficient in some cases, which can enhance the computing efficiency. We fix $k_0=5$, other parameters and relative errors are shown in Table \ref{table:ex2-HDG-n123456}. Compared to the PINNs based on multilevel domain decomposition in \cite{Dolean2023Multilevel}, we see that the LRNN-HDG method can reach more accurate numerical solutions.

\begin{example} [Stokes equations with different viscosities] \label{ex3}
Let $\Omega=(0,~1)^2$ and coefficient of viscosity $\nu=0.1$, $0.01$ and $0.001$, the exact solution be
    \begin{align*}
        \bu(x,~y)=\left(\begin{matrix}\sin^2(\pi x)\sin(2\pi y)\\-\sin(2\pi x)\sin^2(\pi y)\end{matrix}\right), \quad p(x,~y)=\cos(\pi x)\cos(\pi y).
    \end{align*}
\end{example}

We determine the numbers of neurons by the index $k_0$ as follows: 
\begin{align*}
    N_{\underline{\bsigma}_r}=\frac{1}{2}(k_0+1)(k_0+2),\quad N_{\hat{\underline{\bsigma}}_r}=k_0+1,\quad N_{\bu_r}=\frac{1}{2}(k_0+2)(k_0+3).
\end{align*}

\begin{table}[!htbp]
\centering
\resizebox{\textwidth}{!}{
\begin{tabular}{|c|c|c|c|c|c|c|c|c|c|}
\hline
& \multicolumn{3}{c|}{$h=2^{-1}$, $r=1.5$} & \multicolumn{3}{c|}{$h=3^{-1}$, $r=1.2$} & \multicolumn{3}{c|}{$h=4^{-1}$, $r=0.9$} \\ 
\hline
\diagbox{$k_0$}{err} & $e_0(\underline{\bsigma})$ & $e_0(\bu)$ & $e_0(p)$ & $e_0(\underline{\bsigma})$ & $e_0(\bu)$ & $e_0(p)$ & $e_0(\underline{\bsigma})$ & $e_0(\bu)$ & $e_0(p)$ \\
\hline
4 & 3.46e+00 & 1.11e-01 & 6.40e+00 & 3.37e+00 & 1.48e-02 & 7.39e+00 & 7.62e-01 & 3.47e-03 & 1.63e+00 \\ 
\hline                                                                                                  
5 & 1.11e+00 & 2.40e-02 & 2.43e+00 & 1.02e-01 & 1.99e-03 & 2.23e-01 & 7.48e-03 & 2.11e-04 & 1.60e-02 \\ 
\hline                                                                                                  
6 & 7.80e-02 & 3.11e-03 & 1.71e-01 & 1.78e-03 & 1.26e-04 & 3.74e-03 & 1.46e-04 & 1.36e-05 & 3.02e-04 \\ 
\hline                                                                                                  
7 & 2.08e-03 & 3.14e-04 & 4.30e-03 & 1.14e-04 & 1.76e-05 & 2.21e-04 & 2.16e-05 & 3.30e-06 & 3.84e-05 \\ 
\hline                                                                                                  
8 & 4.67e-04 & 9.88e-05 & 8.99e-04 & 4.26e-05 & 7.54e-06 & 7.53e-05 & 1.16e-05 & 1.76e-06 & 1.89e-05 \\ 
\hline                                                                                                  
9 & 2.35e-04 & 5.57e-05 & 4.20e-04 & 2.85e-05 & 4.69e-06 & 4.89e-05 & 8.36e-06 & 1.33e-06 & 1.37e-05 \\ 
\hline                                                                                                  
10 & 2.03e-04 & 5.18e-05 & 3.50e-04 & 2.12e-05 & 3.75e-06 & 3.59e-05 & 7.49e-06 & 1.11e-06 & 1.21e-05 \\
\hline
\end{tabular}
}
\caption{Relative $L^2$ errors of the LRNN-HDPG method for different $k_0$ and $h$ when $\nu=0.1$ in Example \ref{ex3}. }
\label{table:ex3-nu-1}
\end{table}

\begin{table}[!htbp]
\centering
\resizebox{\textwidth}{!}{
\begin{tabular}{|c|c|c|c|c|c|c|c|c|c|}
\hline
& \multicolumn{3}{c|}{$h=2^{-1}$, $r=1.5$} & \multicolumn{3}{c|}{$h=3^{-1}$, $r=1.2$} & \multicolumn{3}{c|}{$h=4^{-1}$, $r=0.9$} \\ 
\hline
\diagbox{$k_0$}{err} & $e_0(\underline{\bsigma})$ & $e_0(\bu)$ & $e_0(p)$ & $e_0(\underline{\bsigma})$ & $e_0(\bu)$ & $e_0(p)$ & $e_0(\underline{\bsigma})$ & $e_0(\bu)$ & $e_0(p)$ \\
\hline
4 & 3.30e-01 & 1.54e-01 & 4.85e-01 & 4.29e-01 & 1.81e-02 & 6.34e-01 & 1.32e-01 & 3.95e-03 & 1.95e-01 \\ 
\hline                                                                                                  
5 & 1.96e-01 & 2.20e-02 & 2.90e-01 & 2.16e-02 & 2.57e-03 & 3.18e-02 & 1.25e-03 & 2.13e-04 & 1.82e-03 \\ 
\hline                                                                                                  
6 & 1.04e-02 & 3.06e-03 & 1.49e-02 & 2.64e-04 & 1.39e-04 & 3.61e-04 & 4.14e-05 & 2.85e-05 & 5.79e-05 \\ 
\hline                                                                                                  
7 & 7.83e-04 & 8.01e-04 & 1.12e-03 & 7.02e-05 & 6.19e-05 & 9.88e-05 & 1.53e-05 & 1.33e-05 & 2.08e-05 \\ 
\hline                                                                                                  
8 & 3.37e-04 & 4.44e-04 & 4.70e-04 & 4.19e-05 & 4.22e-05 & 5.79e-05 & 1.06e-05 & 9.19e-06 & 1.44e-05 \\ 
\hline                                                                                                  
9 & 1.75e-04 & 2.45e-04 & 2.42e-04 & 2.40e-05 & 2.59e-05 & 3.26e-05 & 7.47e-06 & 6.70e-06 & 1.00e-05 \\ 
\hline                                                                                                  
10 & 1.52e-04 & 2.15e-04 & 2.10e-04 & 1.99e-05 & 2.28e-05 & 2.67e-05 & 7.63e-06 & 6.99e-06 & 1.03e-05 \\
\hline
\end{tabular}
}
\caption{Relative $L^2$ errors of the LRNN-HDPG method for different $k_0$ and $h$ when $\nu=0.01$ in Example \ref{ex3}. }
\label{table:ex3-nu-2}
\end{table}

\begin{table}[!htbp]
\centering
\resizebox{\textwidth}{!}{
\begin{tabular}{|c|c|c|c|c|c|c|c|c|c|}
\hline
& \multicolumn{3}{c|}{$h=2^{-1}$, $r=1.5$} & \multicolumn{3}{c|}{$h=3^{-1}$, $r=1.2$} & \multicolumn{3}{c|}{$h=4^{-1}$, $r=0.9$} \\ 
\hline
\diagbox{$k_0$}{err} & $e_0(\underline{\bsigma})$ & $e_0(\bu)$ & $e_0(p)$ & $e_0(\underline{\bsigma})$ & $e_0(\bu)$ & $e_0(p)$ & $e_0(\underline{\bsigma})$ & $e_0(\bu)$ & $e_0(p)$ \\
\hline
4 & 3.27e-02 & 1.40e-01 & 4.49e-02 & 4.23e-02 & 2.12e-02 & 5.79e-02 & 1.35e-02 & 4.01e-03 & 1.91e-02 \\ 
\hline                                                                                                  
5 & 1.54e-02 & 2.48e-02 & 2.12e-02 & 1.18e-03 & 2.03e-03 & 1.67e-03 & 1.21e-04 & 3.74e-04 & 1.69e-04 \\ 
\hline                                                                                                  
6 & 1.27e-03 & 7.25e-03 & 1.76e-03 & 1.21e-04 & 5.77e-04 & 1.70e-04 & 2.71e-05 & 1.55e-04 & 3.76e-05 \\ 
\hline                                                                                                  
7 & 5.23e-04 & 3.79e-03 & 7.36e-04 & 5.78e-05 & 3.78e-04 & 8.07e-05 & 2.04e-05 & 1.14e-04 & 2.83e-05 \\ 
\hline                                                                                                  
8 & 2.85e-04 & 1.98e-03 & 4.01e-04 & 3.68e-05 & 2.61e-04 & 5.11e-05 & 1.61e-05 & 8.88e-05 & 2.23e-05 \\ 
\hline                                                                                                  
9 & 1.92e-04 & 1.80e-03 & 2.68e-04 & 3.42e-05 & 2.31e-04 & 4.76e-05 & 1.66e-05 & 8.87e-05 & 2.31e-05 \\ 
\hline                                                                                                  
10 & 1.52e-04 & 1.37e-03 & 2.13e-04 & 2.89e-05 & 2.00e-04 & 4.01e-05 & 1.37e-05 & 7.67e-05 & 1.89e-05 \\
\hline
\end{tabular}
}
\caption{Relative $L^2$ errors of the LRNN-HDPG method for different $k_0$ and $h$ when $\nu=0.001$ in Example \ref{ex3}. }
\label{table:ex3-nu-3}
\end{table}

We solve this problem on several uniform rectangle partitions. In all subsequent examples, we set the test function degree $k=k_0$. The relative $L^2$ errors of $\underline{\bsigma}$, $\bu_r$ and $p_r$ for different $\nu$ and mesh sizes $h$ are shown in Table \ref{table:ex3-nu-1} -- Table \ref{table:ex3-nu-3}. We see that the LRNN-HDPG method achieves good results for various $\nu=0.1$, $0.01$, and $0.001$.

\begin{example}[Stokes-Darcy equations with Beavers-Joseph law]  \label{ex4}
Let $\underline{\bK}=\alpha\underline{\bI}_{2\times2}$, $\kappa=1$, $\Omega_S=(0,~\pi)\times(0,~\pi)$, $\Omega_D=(0,~\pi)\times(-\pi,~0)$ and $\Omega=\Omega_S\cup\Omega_D$. The exact solution is 
    \begin{align*}
        \bu^S(x,~y)=\alpha\left(\begin{matrix}\sin(2y)\cos(x)\\(\sin^2(y )-2)\sin(x)\end{matrix}\right),\quad p^S(x,~y)=\sin(2x)\sin(2y)
    \end{align*}
    and
    \begin{align*}
        p^D(x,~y)=\frac{2}{\pi}y(\pi-y)\sin(x).
    \end{align*}
    The corresponding $\fb^S$, $\bg^S$, $f^D$ and $g^D$ can be calculated directly. We choose $\alpha=0.01,~100$ and $\nu=0.001,~0.1$. We consider the Stokes-Darcy equations with interface conditions \eqref{Eq:interface1}, \eqref{Eq:interface2} and the Beavers-Joseph law \eqref{Eq:interface5}. 
\end{example}

\begin{table}[!htpb]
\centering
\begin{tabular}{|c||c|c|c|c||c|c|c|c||}
\hline
& \multicolumn{4}{c||}{$\alpha=0.01,~\nu=0.001,~r^D=0.7,~r^S=0.6$} & \multicolumn{4}{c||}{$\alpha=0.01,~\nu=0.1,~r^D=0.7,~r^S=0.6$} \\ 
\hline
\diagbox{$k_0$}{err} & $e_0(\bu^S)$ & $e_0(p^S)$ & $e_0(\bu^D)$ & $e_0(p^D)$ & $e_0(\bu^S)$ & $e_0(p^S)$ & $e_0(\bu^D)$ & $e_0(p^D)$ \\
\hline
4 & 1.27e-01 & 4.03e-03 & 1.67e-03 & 1.63e-04 & 1.39e-02 & 4.02e-03 & 2.00e-03 & 1.11e-04 \\
\hline
5 & 8.60e-02 & 1.94e-03 & 5.14e-04 & 6.69e-05 & 3.17e-03 & 1.05e-03 & 2.87e-04 & 2.93e-05 \\
\hline
6 & 4.11e-02 & 1.32e-03 & 3.25e-04 & 3.74e-05 & 2.02e-03 & 7.57e-04 & 1.70e-04 & 1.86e-05 \\
\hline
7 & 3.36e-02 & 1.03e-03 & 3.18e-04 & 3.62e-05 & 1.83e-03 & 5.50e-04 & 1.42e-04 & 1.51e-05 \\
\hline
8 & 3.17e-02 & 8.17e-04 & 3.23e-04 & 3.62e-05 & 1.86e-03 & 5.21e-04 & 1.16e-04 & 1.22e-05 \\
\hline
9 & 2.61e-02 & 8.03e-04 & 2.83e-04 & 3.23e-05 & 1.80e-03 & 4.62e-04 & 1.21e-04 & 1.26e-05 \\
\hline
\hline
& \multicolumn{4}{c||}{$\alpha=100,~\nu=0.001,~r^D=1.3,~r^S=1.5$} & \multicolumn{4}{c||}{$\alpha=100,~\nu=0.1,~r^D=1.3,~r^S=1.5$} \\ 
\hline
\diagbox{$k_0$}{err} & $e_0(\bu^S)$ & $e_0(p^S)$ & $e_0(\bu^D)$ & $e_0(p^D)$ & $e_0(\bu^S)$ & $e_0(p^S)$ & $e_0(\bu^D)$ & $e_0(p^D)$ \\
\hline
4 & 4.36e-03 & 9.11e-02 & 2.85e-03 & 3.49e-04 & 3.02e-03 & 1.62e+01 & 2.87e-02 & 1.02e-03 \\
\hline
5 & 1.63e-03 & 1.72e-02 & 3.39e-03 & 2.44e-04 & 8.92e-04 & 1.05e+00 & 3.99e-03 & 2.46e-04 \\
\hline
6 & 1.16e-03 & 1.18e-02 & 8.56e-04 & 8.39e-05 & 7.95e-05 & 1.78e-02 & 2.56e-04 & 1.82e-05 \\
\hline
7 & 4.16e-04 & 5.05e-03 & 6.66e-04 & 9.09e-05 & 1.77e-05 & 5.96e-03 & 1.06e-04 & 8.12e-06 \\
\hline
8 & 4.34e-04 & 6.35e-03 & 5.25e-04 & 8.11e-05 & 1.71e-05 & 4.41e-03 & 9.13e-05 & 9.62e-06 \\
\hline
9 & 3.56e-04 & 4.80e-03 & 4.61e-04 & 7.16e-05 & 5.33e-06 & 1.63e-03 & 8.71e-05 & 1.07e-05 \\
\hline
\end{tabular}
\caption{Relative $L^2$ errors of the LRNN-HDPG method with $M=30$ for different $\alpha$ and $\nu$ in Example \ref{ex4}. }
\label{table:ex4-alpha-nu}
\end{table}

For a given positive integer $k_0$, we set
\begin{align*}
    N_{\underline{\bsigma}_r^S}=\frac{1}{2}(k_0+2)(k_0+3),\quad N_{\hat{\underline{\bsigma}}_r^S}=k_0+2,\quad N_{\bu_r^S}=\frac{1}{2}(k_0+3)(k_0+4),\\
    N_{\bu_r^D}=\frac{1}{2}(k_0+1)(k_0+2),\quad N_{\hat{\bu}_r^D}=k_0+1,\quad N_{p_r^D}=\frac{1}{2}(k_0+2)(k_0+3),
\end{align*}
and the test function degrees $k^S=k_0+1$ and $k^D=k_0$. 
We solve this problem on uniform rectangle partitions with $2N_x=N_y=6$. For LRNN spaces defined on different domains, we have different uniform distributions, denoted by $\cU(-r^D,~r^D)$ and $\cU(-r^S,~r^S)$. We recall that $M$ is the number of sampling points on $\Gamma$, and we take $M=30$ in this example. The relative $L^2$ errors of $\bu_r^S$, $p_r^S$, $\bu_r^D$ and $p_r^D$ for different $\alpha$, $\nu$ and $k_0$ are shown in Table \ref{table:ex4-alpha-nu}. 
For different choices of $\alpha$ and $\nu$, the LRNN-HDPG method shows good performance.

\begin{example}[Brinkman equations] \label{ex5}
Let $\Omega=(0,~1)^2$, $\underline{\bkappa}^{-1}=\alpha\underline{\bI}_{2\times2}$ and $\nu=1$, where $\alpha$ is a positive constant. The exact solution is
    \begin{align*}
        \bu(x,~y)=\left(\begin{matrix}x^2(x-1)^2y(y-1)(2y-1)\\-x(x-1)(2x-1)y^2(y-1)^2\end{matrix}\right), \quad p(x,~y)=(2x-1)(2y-1).
    \end{align*}
\end{example}

We determine the numbers of neurons by the index $k_0$ as follows: 
\begin{align*}
    N_{\underline{\bsigma}_r}=\frac{1}{2}(k_0+1)(k_0+2),\quad N_{\hat{\underline{\bsigma}}_r}=k_0+1,\quad N_{\bu_r}=\frac{1}{2}(k_0+2)(k_0+3).
\end{align*}
We set the test function degree $k=k_0$.

We solve this problem on uniform rectangle partitions with different mesh sizes. The results are shown in Table \ref{table:ex5-alpha-1e-3} -- \ref{table:ex5-alpha-1e+3}. For $\alpha=10^{-3},~1$ and $10^3$, we report the relative $L^2$ errors of $\underline{\bsigma}_r$, $\bu_r$ and $p_r$ for LRNN-HDPG method. For different mesh sizes $h$, we have different uniform distributions $\cU(-r,~r)$ with $r=0.8$ for $h=2^{-1}$ and $r=0.9$ for $h=3^{-1},~4^{-1}$. We observe that the LRNN-HDPG method performs well for various $\alpha$.

\begin{table}[!htbp]
\centering
\resizebox{\textwidth}{!}{
\begin{tabular}{|c|c|c|c|c|c|c|c|c|c|}
\hline
& \multicolumn{3}{c|}{$h=2^{-1}$, $r=0.8$} & \multicolumn{3}{c|}{$h=3^{-1}$, $r=0.9$} & \multicolumn{3}{c|}{$h=4^{-1}$, $r=0.9$} \\ 
\hline
\diagbox{$k_0$}{err} & $e_0(\underline{\bsigma})$ & $e_0(\bu)$ & $e_0(p)$ & $e_0(\underline{\bsigma})$ & $e_0(\bu)$ & $e_0(p)$ & $e_0(\underline{\bsigma})$ & $e_0(\bu)$ & $e_0(p)$ \\
\hline
3 & 1.23e+00 & 1.00e+00 & 4.11e-01 & 3.22e+01 & 1.02e-01 & 1.08e+01 & 1.27e+01 & 3.81e-02 & 4.26e+00 \\
\hline                                                                                                 
4 & 3.50e-01 & 1.71e-01 & 1.16e-01 & 2.91e-01 & 1.97e-02 & 9.72e-02 & 5.65e-02 & 4.92e-03 & 1.86e-02 \\
\hline                                                                                                 
5 & 4.73e-02 & 1.75e-02 & 1.57e-02 & 3.72e-03 & 1.34e-03 & 1.21e-03 & 5.56e-04 & 2.25e-04 & 1.81e-04 \\
\hline                                                                                                 
6 & 1.41e-03 & 1.07e-03 & 4.54e-04 & 9.97e-05 & 8.71e-05 & 2.90e-05 & 4.10e-05 & 2.99e-05 & 1.13e-05 \\
\hline                                                                                                 
7 & 1.99e-04 & 2.31e-04 & 5.51e-05 & 3.43e-05 & 3.08e-05 & 8.23e-06 & 1.39e-05 & 1.00e-05 & 3.25e-06 \\
\hline                                                                                                 
8 & 9.84e-05 & 1.27e-04 & 2.25e-05 & 2.16e-05 & 2.00e-05 & 5.06e-06 & 1.07e-05 & 7.40e-06 & 2.49e-06 \\
\hline
\end{tabular}
}
\caption{Relative $L^2$ errors of the LRNN-HDPG method for different $k_0$ and $h$ when $\alpha=10^{-3}$ in Example \ref{ex5}. }
\label{table:ex5-alpha-1e-3}
\end{table}

\begin{table}[!htbp]
\centering
\resizebox{\textwidth}{!}{
\begin{tabular}{|c|c|c|c|c|c|c|c|c|c|}
\hline
& \multicolumn{3}{c|}{$h=2^{-1}$, $r=0.8$} & \multicolumn{3}{c|}{$h=3^{-1}$, $r=0.9$} & \multicolumn{3}{c|}{$h=4^{-1}$, $r=0.9$} \\ 
\hline
\diagbox{$k_0$}{err} & $e_0(\underline{\bsigma})$ & $e_0(\bu)$ & $e_0(p)$ & $e_0(\underline{\bsigma})$ & $e_0(\bu)$ & $e_0(p)$ & $e_0(\underline{\bsigma})$ & $e_0(\bu)$ & $e_0(p)$ \\
\hline
3 & 2.36e+00 & 1.00e+00 & 2.36e+00 & 1.61e+01 & 1.05e-01 & 1.61e+01 & 4.65e+01 & 4.89e-02 & 4.66e+01 \\
\hline                                                                                                 
4 & 5.40e-01 & 1.54e-01 & 5.35e-01 & 2.38e-01 & 1.97e-02 & 2.28e-01 & 9.32e-02 & 6.07e-03 & 9.34e-02 \\
\hline                                                                                                 
5 & 7.34e-02 & 1.74e-02 & 7.20e-02 & 4.53e-03 & 1.46e-03 & 4.44e-03 & 9.06e-04 & 3.03e-04 & 8.72e-04 \\
\hline                                                                                                 
6 & 1.52e-03 & 1.43e-03 & 1.45e-03 & 1.19e-04 & 9.37e-05 & 1.08e-04 & 4.27e-05 & 2.74e-05 & 3.66e-05 \\
\hline                                                                                                 
7 & 2.02e-04 & 2.24e-04 & 1.61e-04 & 3.60e-05 & 3.17e-05 & 2.81e-05 & 1.39e-05 & 9.79e-06 & 1.02e-05 \\
\hline                                                                                                 
8 & 1.02e-04 & 1.25e-04 & 7.01e-05 & 2.22e-05 & 2.20e-05 & 1.55e-05 & 9.48e-06 & 6.61e-06 & 6.76e-06 \\
\hline
\end{tabular}
}
\caption{Relative $L^2$ errors of the LRNN-HDPG method for different $k_0$ and $h$ when $\alpha=1$ in Example \ref{ex5}. }
\label{table:ex5-alpha-1e0}
\end{table}

\begin{table}[!htbp]
\centering
\resizebox{\textwidth}{!}{
\begin{tabular}{|c|c|c|c|c|c|c|c|c|c|}
\hline
& \multicolumn{3}{c|}{$h=2^{-1}$, $r=0.8$} & \multicolumn{3}{c|}{$h=3^{-1}$, $r=0.9$} & \multicolumn{3}{c|}{$h=4^{-1}$, $r=0.9$} \\ 
\hline
\diagbox{$k_0$}{err} & $e_0(\underline{\bsigma})$ & $e_0(\bu)$ & $e_0(p)$ & $e_0(\underline{\bsigma})$ & $e_0(\bu)$ & $e_0(p)$ & $e_0(\underline{\bsigma})$ & $e_0(\bu)$ & $e_0(p)$ \\
\hline
3 & 2.35e+01 & 1.00e+00 & 2.36e+01 & 1.96e+00 & 1.08e-01 & 1.96e+00 & 1.28e+00 & 3.73e-02 & 1.28e+00 \\
\hline                                                                                                 
4 & 7.16e-01 & 8.63e-02 & 7.18e-01 & 4.75e-02 & 8.69e-03 & 4.75e-02 & 8.67e-03 & 1.95e-03 & 8.58e-03 \\
\hline                                                                                                 
5 & 1.42e-02 & 6.72e-03 & 1.40e-02 & 1.42e-03 & 7.68e-04 & 1.36e-03 & 3.02e-04 & 1.60e-04 & 2.75e-04 \\
\hline                                                                                                 
6 & 2.20e-03 & 1.77e-03 & 1.98e-03 & 2.77e-04 & 2.29e-04 & 2.27e-04 & 1.08e-04 & 7.32e-05 & 8.35e-05 \\
\hline                                                                                                 
7 & 7.67e-04 & 7.59e-04 & 5.72e-04 & 1.63e-04 & 1.55e-04 & 1.21e-04 & 7.01e-05 & 5.59e-05 & 5.21e-05 \\
\hline                                                                                                 
8 & 7.52e-04 & 7.54e-04 & 5.61e-04 & 1.53e-04 & 1.46e-04 & 1.14e-04 & 6.52e-05 & 5.88e-05 & 4.94e-05 \\
\hline
\end{tabular}
}
\caption{Relative $L^2$ errors of the LRNN-HDPG method for different $k_0$ and $h$ when $\alpha=10^3$ in Example \ref{ex5}. }
\label{table:ex5-alpha-1e+3}
\end{table}

\section{Summary}

We propose the local randomized neural networks with the hybridized discontinuous Petrov-Galerkin methods for solving Darcy flows, Stokes-Darcy problems, and Brinkman equations. The LRNN-HDPG methods are stable in both the Stokes and the Darcy regions, and they can directly enforce the interface conditions at sampling points on the interface. Various numerical examples show that the methods are stable for different coefficients, and can achieve high accuracy with a few degrees of freedom. However, these methods need more research in several aspects. How can we apply the methods to nonlinear partial differential equations? How can we improve the neural network structures to increase the approximation abilities? How can we use parallel processing to enhance their efficiency? How can we use mesh adaptation to improve their performance for complex problems? Moreover, deriving error estimates for the methods is another important topic for future work.


\end{document}